\documentclass{amsart}
\usepackage{amssymb}
\usepackage{graphics}
\usepackage{amscd}
\newtheorem{thm}{Theorem}[section]
\newtheorem{cor}[thm]{Corollary}
\newtheorem{lem}[thm]{Lemma}
\newtheorem{prop}[thm]{Proposition}
\theoremstyle{definition}
\newtheorem{defn}[thm]{Definition}

\theoremstyle{remark}
\newtheorem{rem}[thm]{Remark}
\numberwithin{equation}{section}
\newtheorem{note}[thm]{Notation}
\newcommand{\norm}[1]{\left\Vert#1\right\Vert}
\newcommand{\abs}[1]{\left\vert#1\right\vert}

\newcommand{\Real}{\mathbb R}

\newcommand{\To}{\longrightarrow}

\newcommand{\im}{\mathop{\textrm{Im}}\nolimits}
\newcommand{\imm}{\mathop{\textrm{Imm}}\nolimits}
\begin{document}

\title{Regular homotopy classes of singular maps}%
\author{Andr\'as Juh\'asz}%
\address{Department of Analysis, E\"otv\"os Lor\'and University,
P\'azm\'any P\'eter s\'et\'any 1/C, Budapest, Hungary 1117}%
\email{juhasz.6@dpg.hu}%

\thanks{Research partially supported by OTKA grant no. T037735}%
\subjclass{57R45; 58K30; 57R42}%
\keywords{Whitney umbrella, immersion, regular homotopy, obstruction}%

\date{\today}%
%\dedicatory{}%
%\commby{}%
% ----------------------------------------------------------------
\begin{abstract}
Two locally generic maps $f,g \colon M^n \To \Real^{2n-1}$ are
regularly homotopic if they lie in the same path-component of the
space of locally generic maps. Our main result is that if $n \neq
3$ and $M^n$ is a closed $n$-manifold then the regular homotopy
class of every locally generic map $f \colon M^n \To \Real^{2n-1}$
is completely determined by the number of its singular points
provided that $f$ is \emph{singular} (i.e., $f$ is not an
immersion). This extends the analogous result of \cite{Juhasz} for
$n=2$.
\end{abstract}

\maketitle
% ----------------------------------------------------------------
\section{Introduction}

\begin{figure}[t]
\includegraphics{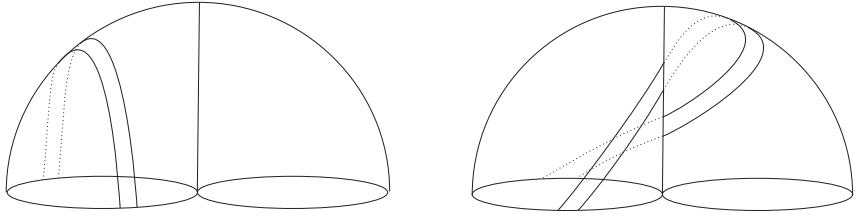}
\caption{Pushing a $1$-cell through a cross-cap point}
\label{fig:1}
\end{figure}

This paper is a sequel to \cite{Juhasz} in which locally generic
maps of closed surfaces into $\Real^3$ were classified up to
regular homotopy. It turned out that for maps with at least one
cross-cap point the number of singular points was the only regular
homotopy invariant. The obstruction to constructing a regular
homotopy was destroyed by pushing $1$-cells of the surface through
a singular point (see Figure \ref{fig:1}).

In the present work we extend this result for locally generic maps
of closed $n$-manifolds into $\Real^{2n-1}$ in case $n \neq 3$
(Theorem \ref{thm:1}). However, this is not simply an adaptation
of the ideas of \cite{Juhasz}. The idea of pushing
$(n-1)$-simplices through a singular point originates from the
previous work, but a lot of technical problems have to be dealt
with. For example, if $n$ is odd the Whitney-umbrella points have
signs and the $(n-1)$-dimensional obstruction is
$\mathbb{Z}$-valued. Moreover, a new type of obstruction appears
if $n > 3$ that was not present in the case $n = 2$. This is an
$n$-dimensional obstruction and is related to the double point set
of the map. In eliminating this obstruction we make essential use
of a result of T. Ekholm \cite{Ekholm2}.

It seems that this second type of obstruction cannot be eliminated
if $n=3$. In fact, in addition to the number of singular points a
new $\mathbb{Z}$-valued invariant comes into the picture. I intend
to deal with this in a separate paper.

\subsection{Preliminary definitions}

\begin{defn}
Let $M^n$ be a closed $n$-manifold and $N^{2n-1}$ an arbitrary
$(2n-1)$-manifold . A map $f \colon M^n \to N^{2n-1}$ is called
\emph{locally generic} if it is an immersion except for cross-cap
(or Whitney-umbrella) singularities. Thus for any singular point
$p$ there exist coordinate systems $(x_1, \dots, x_n)$ about $p$
and $(y_1, \dots, y_{2n-1})$ about $f(p)$ such that $f$ is given
near $p$ by \begin{equation} \label{eqn:1} y_1 = x_1^2, \,\, y_i =
x_i, \,\, y_{n+i-1} = x_1 x_i \,\, (i = 2, \dots, n).
\end{equation}

\end{defn}

\begin{note}
Let $L(M^n, N^{2n-1})$ denote the subspace of locally generic
maps, $G(M^n, N^{2n-1})$ the subspace of generic maps and
$\imm(M^n, N^{2n-1})$ the subspace of immersions in
$C^{\infty}(M^n, N^{2n-1})$ endowed with the $C^{\infty}$
topology. We introduce the notation $L'(M^n, N^{2n-1})$ for the
space $$L(M^n, N^{2n-1}) \setminus \imm(M^n, N^{2n-1}).$$ For $f
\in L(M^n, N^{2n-1})$ the set of singular points of the map $f$ is
denoted by $S(f)$.
\end{note}

\begin{rem} \label{rem:1}
Since $M^n$ is compact and the singularities of $f$ are isolated
we have $\abs{S(f)}<\infty$. It is well known that for $M^n$
closed $\abs{S(f)}$ is an even number.

For $n>2$ a map $f \colon M^n \to N^{2n-1}$ is generic iff
\begin{itemize}
\item it is an immersion with normal crossings except in a finite
set of points
\item the singular points are non-multiple cross-cap points
\item $f$ has at most double crossings.
\end{itemize}

For $n=2$ a generic $f$ might also have triple points.
\end{rem}

\begin{defn}
A \emph{regular homotopy} is a path in the space $L(M^n,
N^{2n-1})$. In other words, it is a smooth map $H \colon M^n
\times [0,1] \to N^{2n-1}$ such that $H_t \in L(M^n, N^{2n-1})$
for every $t \in [0,1]$. Here $H_t(x) = H(x,t)$ for $x \in M^n$
and $t \in [0,1]$. The regular homotopy $H$ is called
\emph{singularity fixing} if $S(H_t) = S(H_0)$ for every $t \in
[0,1]$.
\end{defn}

\begin{defn}
Two maps $f,g \in L(M^n, N^{2n-1})$ are called \emph{regularly
homotopic} if there exists a regular homotopy $H$ such that $H_0 =
f$ and $H_1 = g$. If $H$ can be chosen to be singularity fixing
then we say that $f$ and $g$ are \emph{regularly homotopic through
a singularity fixing homotopy}.
\end{defn}

\begin{note}
The fact that $f$ and $g$ are regularly homotopic is denoted by $f
\sim g$. Furthermore, $f \sim_s g$ denotes that $f$ and $g$ are
regularly homotopic through a singularity fixing homotopy.
\end{note}

\subsection{The classification theorem}

\begin{thm} \label{thm:1}
Suppose that $n > 3$ or $n=2$. Let $M^n$ be a closed $n$-manifold
and $f,g \in L'(M^n, \Real^{2n-1})$. Then $$f \sim g
\Leftrightarrow \abs{S(f)} = \abs{S(g)}.$$
\end{thm}

\section{Known results used in the proof of Theorem \ref{thm:1}}

\subsection{Smale's lemma and $M^n$-regular homotopies}

We will use extensively the following result of Smale (this is
Theorem 1.1 in \cite{Hirsch}, the original paper of Smale is
\cite{Smale}). First we need a few definitions.

\begin{defn} \label{defn:1}
Let $\mathcal{E} = \mathcal{E}_{k,m}$ be the space of all
$C^{\infty}$ immersions of $D^k$ in $\Real^m$ in the $C^1$
topology. Let $\mathcal{B} = \mathcal{B}_{k,m}$ be the set of
pairs $(g, g')$ where $g \colon S^{k-1} \to \Real^m$ is a
$C^{\infty}$ immersion and $g' \colon S^{k-1} \to T \Real^m$ is a
$C^{\infty}$ transversal field of $g$. $\mathcal{B}$ is
topologized as a subspace of Cartesian product of the space of
immersions $S^{k-1} \to \Real^m$, in the $C^1$ topology, with the
space of continuous maps $S^{k-1} \to T\Real^m$, in the
compact-open topology.

If $h \in \mathcal{E}$ then let $h' \colon S^{k-1} \to T\Real^m$
be defined by $h'(x) =$ derivative of $h$ along the radius at $x
\in S^{k-1}$. I.e., if $r(x)$ is the unit tangent vector of $D$
that is normal to $S^{k-1}$ at $x$ and which points away from the
origin then $h'(x) = h_*r(x)$. We define the map $\pi \colon
\mathcal{E} \to \mathcal{B}$ by the formula $\pi(h) = (h|S^{k-1},
h')$. It is clear that $\pi$ is continuous. The following theorem
is Theorem 1.1 in \cite{Hirsch}.
\end{defn}

\begin{thm} \label{thm:9}
If $k<m$, then $\pi \colon \mathcal{E}_{k,m} \to
\mathcal{B}_{k,m}$ has the covering homotopy property.
\end{thm}

The intuitive content of this theorem is as follows: If we are
given an immersed disk $D^k$ in $\Real^m$ such that $k < m$ and we
deform the boundary of the disk and the normal derivatives along
the boundary, then we can deform the whole disk at the same time
so as to induce the given deformation on the boundary and normal
derivatives.

The following definitions were introduced by Hirsch in
\cite{Hirsch}.

\begin{defn}
Let $A$ be an arbitrary subset of the manifold $M^n$. Let $h
\colon A \to Q^q$ and $h' \colon TM^n|A \to TQ^q$ ($Q^q$ is a
manifold) be continuous maps such that $h'$ covers $h$. The pair
$(h, h')$ is called an \emph{$M^n$-regular map}, or
\emph{$M^n$-immersion}, \emph{of $A$ in $Q^q$} if the following
condition is satisfied: there is a neighborhood $V$ of $A$ in
$M^n$ and an immersion $l \colon V \to Q^q$ such that $dl|(TM|A) =
h'$. It follows that $l|A = h$. We say that $(h, h')$ is $C^k$ if
$l$ can be chosen to be $C^k$.
\end{defn}

\begin{defn}
Let $B \subset A \subset M^n$ be subsets. If $(r,r'), (s,s')
\colon A \to Q^q$ are $M^n$-immersions such that $r|B = s|B$ and
$r'|(TM^n|B) = s'|(TM^n|B)$, we say that \emph{$(r,r')$ and
$(s,s')$ are tangent on $B$}, and write this as $(r,r')|B =
(s,s')|B$.
\end{defn}

\begin{defn}
Let $(r,r')$ and $(s,s')$ be $M^n$-immersions of $A$ in $Q^q$ such
that for a certain (possibly empty) subset $B$ of $A$, $(r,r')$
and $(s,s')$ are tangent on $B$. We say that $(r,r')$ and $(s,s')$
are \emph{$M^n$-regularly homotopic (rel B)} if there is a path
$(h_t, h'_t)$ in the space of all $M^n$-immersions of $A$ in $Q^q$
joining $(r,r')$ to $(s,s')$, such that for each $t$, $(h_t,
h_t')|B = (r,r')|B$. Such a path is called an \emph{$M^n$-regular
homotopy (rel B)}, and it is $C^k$ if every $(h_t, h'_t)$ is
$C^k$.
\end{defn}

\begin{note}
The space of all $C^{\infty}$ $\Real^n$-immersions of $D^k$ in
$\Real^q$ is denoted by $\mathcal{I}(k, q; n)$; the space of all
$C^{\infty}$ $\Real^n$-immersions of $S^{k-1}$ in $\Real^q$ is
denoted by $\mathcal{I}'(k, q; n)$. Now let $\pi_n \colon
\mathcal{I}(k, q; n) \to \mathcal{I}'(k, q; n)$ be defined by $$
\pi_n(f,f') = (f|S^{k-1}, f'|(T\Real^n|S^{k-1})).$$ If $n=k$, this
is the map $\pi \colon \mathcal{E}_{k,q} \to \mathcal{B}_{k,q}$
defined in Definition \ref{defn:1}.

For $(f,f') \in \mathcal{I}(k,q;n)$ put
$\Gamma_n(f,f')=\pi_n^{-1}(\pi_n(f,f'))$.
\end{note}

The following statement is Theorem 3.5 in \cite{Hirsch}. It is a
generalization of Smale's lemma to $M^n$-immersions. To avoid
confusion we will refer to it also as Smale's lemma (despite the
fact that it was proved by Hirsch).

\begin{thm} \label{thm:2}
$\pi_n$ has the covering homotopy property if $k < q$.
\end{thm}

Theorem 3.2 in \cite{Hirsch} gives an alternative description of
$M^n$-immersions:

\begin{note}
$ \left\{\, e_1(x), \dots, e_n(x) \,\right\}$ denotes the standard
basis of $T_x\Real^n$.
\end{note}

\begin{lem} \label{lem:1}
There is a homeomorphism $\zeta$ between the space of
$\Real^{n}$-immersions $(h,h') \colon D^k \to \Real^q$ and the
space of pairs $(l, \Psi)$, where $l \in \imm(D^k, \Real^q)$ and
$\Psi$ is a transversal $(n-k)$-field along $l$. The homeomorphism
is given by $\zeta(h,h') = (h, \Psi)$ where $\Psi(x) = h'\{
e_{k+1}(x), \dots, e_n(x) \}$. Moreover, $(f, f')$ is $C^k$ if and
only if $f$ and $\Psi$ are $C^k$.
\end{lem}

An analogous result holds for $S^{k-1}$ (see \cite{Hirsch},
Theorem 3.3):

\begin{note}
Denote the space of pairs $(l, \Psi)$, where $l \in \imm(S^{k-1},
\Real^q)$ and $\Psi$ is a transversal $(n-k+1)$-field along $l$ by
$$\imm_{n-k+1}(S^{k-1}, \Real^q).$$
\end{note}

\begin{lem} \label{lem:2}
There is a homeomorphism  $$\chi \colon \mathcal{I}'(k,q;n) \to
\imm_{n-k+1}(S^{k-1}, \Real^q).$$ $\chi$ is given as follows: Let
$\Phi$ be the normal $(n-k+1)$-field on $S^{k-1}$ given by
$$\Phi(x) = \{r(x), e_{k+1}(x), \dots, e_n(x)\},$$ where $r(x)$ is
the outward  unit normal to $S^{k-1}$ in $\Real^k$. Then
$\chi(h,h') = (h, \Psi)$, where $\Psi(x) = f'(\Phi(x))$.
\end{lem}

\subsection{The obstructions $\tau$ and $\Omega$}

Hirsch defined in \cite{Hirsch} an invariant $\tau(g')$ for each
$(g,g') \in \mathcal{I}'(k,q;n)$. The vanishing of $\tau(g')$
implies that $(g,g')$ comes from $\mathcal{I}(k,q;n)$.

\begin{defn}
$(f,f') \in \mathcal{I}'(k,q;n)$ is said to be \emph{extendible}
if there is a $(g,g') \in \mathcal{I}(k,q;n)$ such that
$\pi_q(g,g') = (f,f')$.
\end{defn}

\begin{defn} \label{defn:3}
Let $(f,f') \colon S^{k-1} \to \Real^n$ be a $C^{\infty}$
$\Real^n$-immersion, i.e., $(f,f') \in \mathcal{I}'(k,q;n)$. The
\emph{obstruction to extending $(f,f')$}, denoted by $\tau(f') \in
\pi_{k-1}(V_{q,n})$ is the homotopy class of the map $S^{k-1} \to
V_{q,n}$ defined by $$x \mapsto f'\{e_1(x), \dots, e_n(x)\}.$$
\end{defn}

The following lemma is Theorem 3.9 in \cite{Hirsch}:

\begin{lem} \label{lem:12}
If $k < q$ and $\tau(f') = 0$ then $(f,f')$ is extendible.
\end{lem}

\begin{defn}
Let $\Phi_n \colon \mathcal{I}(k,q;n) \to C^0(D^k, V_{q,n})$ be as
follows: $$\Phi_n(f,f')(x) = f'\{e_1(x), \dots, e_n(x)\}.$$
\end{defn}

Let $(f,f'), (g,g') \in \mathcal{I}(k,q;n)$, with $\pi_n(f,f') =
\pi_n(g,g')$, so that $(f,f') \in \Gamma_n(g,g')$. Then
$\Phi_n(f,f')$ and $\Phi_n(g,g')$ are maps $D^k \to V_{q,n}$ which
are tangent on $S^{k-1}$.

\begin{defn}
Let $A$ be a topological space, simple in dimension $k$. Let $f,g
\colon D^k \to A$ and assume that $f(x) = g(x)$ if $x \in
S^{k-1}$. Then $d(f,g) \in \pi_k(A)$ is represented by mapping the
"top" hemisphere of $S^k$ by $f$ and the "bottom" one by $g$,
assuming that the orientation of $S^k$ is given by the coordinate
frame $\{\,e_1, \dots, e_k\,\}$ at the "North" pole of $S^k$.
\end{defn}

\begin{defn}
$$\Omega(f',g') = d(\Phi_n(f,f'), \Phi_n(g,g')) \in
\pi_k(V_{q,n})$$is called the \emph{obstruction to an
$\Real^n$-regular homotopy} (rel $S^{k-1}$) between $(f,f')$ and
$(g,g')$. (This is well defined since $V_{q,n}$ is simple in all
dimensions.)
\end{defn}

\begin{rem}
An explicit definition of $\Omega(f',g')$ is as follows: identify
the upper and lower hemispheres of $S^k$ with $D^k$. Let $\omega
\colon S^k \to V_{q,n}$ be the map $\omega(x) = f'\{e_1(x), \dots,
e_n(x) \}$ if $x$ is in the upper hemisphere, $\omega(x) =
g'\{e_1(x), \dots, e_n(x)\}$ if $x$ is in the lower hemisphere.
$\omega(x)$ is well defined on the equator because $(f,f')$ and
$(g,g')$ agree on $S^{k-1}$. Then $\Omega(f',g')$ is the homotopy
class of $\omega$.
\end{rem}

The following lemma is Theorem 4.3 in \cite{Hirsch}:

\begin{lem} \label{lem:9}
Let $(f,f'), (g,g'), (h,h')$ be $C^{\infty}$ $\Real^n$-immersions
of $\Delta^k$ in $\Real^q$ which are all tangent on $\partial
\Delta^k$. Then

(a) $\Omega(f',g')+ \Omega(g',h') = \Omega(f',h')$.

(b) $\Omega(f',f') = 0.$

(c) Given $\alpha \in \pi_k(V_{q,n})$ there exists $(g,g')$ such
that $\Omega(f',g') = \alpha $.

(d) Suppose that $\Omega(f',g') = 0$ and $k<q$. Let $F \colon
\Delta^k \times I \to V_{q,n}$ be a homotopy (rel $\partial
\Delta^k$) between the maps $F_0, F_1 \colon \Delta^k \to V_{q,n}$
defined respectively by $x \mapsto f'\{e_i(x)\}$ and $x \mapsto
g'\{e_i(x)\}$, $i = 1 \dots, n$. Then there is a $C^{\infty}$
$\Real^n$-regular homotopy $(f_t,f_t')$ between $(f,f')$ and
$(g,g')$ such that the map $\Delta^k \times I \to V_{q,n}$ defined
by $(x,t) \mapsto f_t'(e_i(x))$ is homotopic (rel $(\partial
\Delta^k \times I) \cup (\Delta^k \times \partial I)$) to $F$.

\end{lem}

From the proof of Theorem 3.9 in \cite{Hirsch} we can deduce the
following

\begin{thm} \label{thm:3}
Suppose that the $C^{\infty}$ $\Real^n$-immersions $(f,f'), (g,g')
\colon \partial\Delta^n \to \Real^q$ are tangent on
$(\partial\Delta^n) \setminus \text{int}(\Delta^{n-1})$, where
$\Delta^{n-1}$ is an $(n-1)$-face of the $n$-simplex $\Delta^n$.
Then
$$\Omega(f'|\Delta^{n-1}, g'|\Delta^{n-1}) = \tau(f')-\tau(g').
$$
\end{thm}

\subsection{Hirsch's theorem and lemma}

In \cite{Hirsch} Hirsch proves the following theorem that reduces
regular homotopy classification of immersions to homotopy theory.

\begin{thm} \label{thm:4}
Suppose that $n \le q$ and let $M^n$ and $Q^q$ be manifolds such
that $M^n$ is open if $n = q$. Then the natural map $$ \imm(M^n,
Q^q) \to \text{Mono}(TM, TQ) $$ is a weak homotopy equivalence (so
it induces a bijection in $\pi_0$).
\end{thm}

The relative version of Hirsch's theorem also holds:

\begin{thm} \label{thm:5}
Let $M^n$ and $Q^q$ be as in Theorem \ref{thm:4}. Suppose that $K
\subset M$ is compact and $V_K$ is a neighborhood of $K$. Given an
immersion $h \in \imm(V_K, Q^q)$ with differential $$dh = \Phi
\colon TM|V_K \to TQ|V_K,$$ there exists an immersion $l \in
\imm(M^n, Q^q)$ such that $l|K = f$ and $dl$ is homotopic to
$\Phi$.
\end{thm}

As a consequence of Theorem \ref{thm:4} we obtain Hirsch's lemma:

\begin{lem} \label{lem:3}
Let $k \ge 1$ and $f \in \imm(M^n, \Real^{n+k+1})$ be an immersion
with normal bundle $\mu^k \oplus \varepsilon^1$, where
$\varepsilon^1$ denotes the trivial line bundle over $M^n$. Let
$i$ denote the standard embedding $\Real^{n+k} \to \Real^{n+k+1}$.
Then there exists an immersion $g \in \imm(M^n, \Real^{n+k})$ with
normal bundle $\mu^k$ such that $f \sim i \circ g$. During the
regular homotopy a trivialization of the $\varepsilon^1$ component
of $\nu(f \subset \Real^{n+k+1})$ can be deformed simultaneously
in the space of normal fields to be finally vertical along $i
\circ g$.
\end{lem}

\begin{rem}
The relative version of Hirsch's lemma also holds, as can be seen
by using the relative version of Hirsch's theorem.
\end{rem}

\begin{cor} \label{cor:3}
If $k \ge 2$ then Lemma \ref{lem:3} gives a bijection
$$b \colon \pi_0(\imm_1(M^n,
 \Real^{n+k+1})) \to \pi_0(\imm(M^n, \Real^{n+k})). $$
\end{cor}

\subsection{Whitney-umbrellas}

In this chapter we sum up some ideas from \cite{Whitney}. Fix an
orientation for $\Real^{2n-1}$.

\begin{defn}
Suppose that $N^n$ is an $n$-manifold with boundary and that $f
\colon N^n \to \Real^{2n-1}$ is a generic map (thus $f|\partial
N^n$ is an embedding). Let $\nu \colon \partial N^n \to
T\Real^{2n-1}$ be a transversal vector field along $f|\partial
N^n$ such that $-\nu$ points into $f(N^n)$. Then for $\varepsilon$
sufficiently small define $\mathcal{L}_f(N)$ to be the linking
number $$\text{lk}(f|(\partial N^n),f|(\partial N^n) + \epsilon
\nu).$$
\end{defn}

\begin{rem}
Notice that $(f,\nu)$ corresponds to the $N^n$-immersion
$(f,df)|\partial N^n$.
\end{rem}

As a consequence of Lemma 6 in \cite{Whitney} we obtain the
following lemma.

\begin{lem} \label{lem:4}
Suppose that $f \in L(M^n, \Real^{2n-1})$ and $p \in S(f)$. Choose
a coordinate neighborhood $U_p$ of $p$ diffeomorphic to $D^n$ such
that $f|U$ is generic and has the form (\ref{eqn:1}). Then
$\mathcal{L}_f(U_p) = \pm 1$, moreover, the sign does not depend
on the choice of $U_p$ if $n$ is odd. (Recall that for $n$ even
$\mathcal{L}_f$ is defined mod $2$.)
\end{lem}

\begin{note} \label{note:1}
Denote the map $(f,df)|\partial U_p \in \mathcal{I}'(n,2n-1;n)$ by
$(w,w')$.
\end{note}

\begin{cor}
For $n$ even $\mathcal{L}_f(U_p) \equiv 1$ (mod $2$) for every
neighborhood $U_p$.
\end{cor}

In \cite{Whitney} signs are defined for Whitney-umbrellas if $n$
is odd and $\Real^{2n-1}$ is oriented. This can be done as follows
(this is not the original definition of Whitney):

\begin{defn} \label{defn:4}
Suppose that $f \in L(M^n, \Real^{2n-1})$ and $p \in S(f)$. Choose
$U_p$ as in Lemma \ref{lem:4}. The sign of the Whitney-umbrella at
$p$ is then defined as
$$\text{sgn}(p) = \mathcal{L}_f(U_p).$$
\end{defn}

\begin{note}
If $n$ is odd and $f \in L(M^n, \Real^{2n-1})$ denote the set of
positive (respectively negative) cross-cap points of $f$ by
$S_+(f)$ (respectively $S_-(f)$). Moreover let
\[
\#S(f) = \left\{
\begin{array}{ll}
\abs{S_+(f)} - \abs{S_-(f)} & \mbox{if $n$ is odd,}\\
\abs{S(f)} \, \text{mod}\, 2 & \mbox{if $n$ is even,}
\end{array}
\right.
\]
where $\abs{A}$ denotes the cardinality of the set $A$.
\end{note}

In \cite{Whitney} Theorem 2 states the following:

\begin{lem} \label{lem:5}
Let $N^n$ be a compact $n$-manifold with boundary and $$f \in
G(N^n, \Real^{2n-1}).$$ Then $\mathcal{L}_f(N^n) = \#S(f)$.
\end{lem}

\begin{rem}
If $f$ is only locally generic, we may perturb it slightly to
obtain a generic map. From this we can generalize Lemma
\ref{lem:5} to locally generic maps.
\end{rem}

As a trivial consequence we obtain

\begin{cor} \label{cor:1}
Let $M^n$ be a closed $n$-manifold and $f \in L(M^n,
\Real^{2n-1})$. Then for $n$ odd the algebraic number of singular
points vanishes (i.e., $\abs{S_+(f)} = \abs{S_-(f)}$), while for
$n$ even, it vanishes mod $2$ (i.e., $\abs{S(f)}$ is even).
\end{cor}

We will also use Theorem $7$ of \cite{Whitney}, which states the
following:

\begin{thm} \label{thm:8}
Let $M^n$ be a compact and connected manifold with boundary and $f
\colon M^n \to \Real^{2n-1}$ be a smooth map which is an embedding
in a neighborhood of $\partial M^n$. Suppose that
$\mathcal{L}_f(M^n) = 0$ if $n$ is odd and $\mathcal{L}_f(M^n)
\equiv 0$ mod $2$ if $n$ is even. Then there exists an immersion
$g \colon M^n \to \Real^{2n-1}$ that is arbitrarily $C^0$-close to
$f$ and equals $f$ in a neighborhood of $\partial M^n$.
\end{thm}

Finally, we present an equivalent definition of the sign of a
Whitney-umbrella singularity (motivated by \cite{Saeki}). With an
additional choice we can also define the sign of a singular point
if $n$ is even.

\begin{defn} \label{defn:5}
Let $f \in L(M^n, \Real^{2n-1})$ and $p \in S(f)$. Fix an
orientation of $\Real^{2n-1}$. Then we define the sign of $p$ as
follows.

Choose local coordinates $(x_1, \dots, x_n)$ about $p$ and $(y_1,
\dots, y_{2n-1})$ about $f(p)$ such that $f$ is given near $p$ by
equation (\ref{eqn:1}). Let $D_{\varepsilon} = \{\,y_1^2 + \dots +
y_{2n-1}^2 \le \varepsilon \,\}$ for $\varepsilon > 0$
sufficiently small. Then Lemma 2.2 in \cite{Saeki} states that
$B_{\varepsilon} = f^{-1}(D_{\varepsilon})$ is a closed disc
neighborhood of $p$ in $M^n$. Choose an arbitrary orientation of
$S_{\varepsilon} = \partial B_{\varepsilon}$. Denote by $q$ the
double point of $f|S_{\varepsilon}$ and let $f^{-1}(q) = \{\,q_1,
q_2\,\}$. Moreover, let $\underline{v} = (v_1, \dots, v_{n-1})$ be
a positive basis of $T_{q_1}S_{\varepsilon}$ and $\underline{w} =
(w_1, \dots, w_{n-1})$ a positive basis of
$T_{q_2}S_{\varepsilon}$. Orient $\partial D_{\varepsilon}$ by its
outward normal vector. Then $p$ is called positive if
$(df(\underline{v}), df(\underline{w}))$ is a positive basis of
$T_q(\partial D_{\varepsilon})$, and negative otherwise.

This definition does not depend on the choice of the orientation
of $S_{\varepsilon}$. However, if $n$ is even then the sign of $p$
does depend on the ordering of the points $q_1$ and $q_2$, i.e.,
if we swap the two points then the sign of $p$ also changes. For
$n$ odd the sign of $p$ is independent of the ordering of $q_1$
and $q_2$. Thus fixing the orientation of $\Real^{2n-1}$ defines
the sign of a Whitney-umbrella point only if $n$ is odd. If we
also fix an ordering of the two branches of $f$ meeting at the
open double curve ending at $f(p)$ then this defines a sign of $p$
if $n$ is even.

It is easy to verify that Definition \ref{defn:4} and Definition
\ref{defn:5} are equivalent.
\end{defn}

\section{Proof of Theorem \ref{thm:1}}

\subsection{Outline}

The proof of the case $n=2$ can be found in \cite{Juhasz}. The
generalizations of those ideas are incorporated into the present
proof for the case $n>3$.

The implication $$f \sim g \Rightarrow \abs{S(f)} = \abs{S(g)}$$
is Proposition 2.4 in \cite{Juhasz}. The proof relies on the fact
that the cross-cap singularity is stable and $M^n$ is closed. Thus
we will only prove the other implication.

For $n>3$ the proof is divided into two main parts. Using
Hirsch-Smale theory (for a reference see \cite{Hirsch}) we will
see that there are two obstructions to constructing a regular
homotopy between two immersions from $M^n$ to $\Real^{2n-1}$.
These obstructions can be eliminated in the presence of cross-cap
points.

In the first part we will construct a regular homotopy that pushes
the $(n-1)$-simplices of $M^n$ through a singular point to destroy
the first obstruction. In the second part we will first reduce the
problem to the case $M^n = S^n$. To remove the second obstruction
we will merge the closed double curves of $f$ and $g$ with the
ones connecting the singular points. This can be done using a
variant of the Whitney-trick. Finally, we replace $f$ and $g$ with
immersions and we use an argument of \cite{Ekholm2} that shows
that the Smale-invariant of an immersion of $S^n \to \Real^{2n-1}$
is completely determined by the geometry of its self-intersection
if $n \ge 4 $.

\subsection{Setup}

From now on $f$ and $g$ denote the maps in the statement of
Theorem \ref{thm:1}. Thus $f, g \in L'(M^n, \Real^{2n-1})$ and
$\abs{S(f)} = \abs{S(g)}$. We also suppose that $n>3$ and
$\Real^{2n-1}$ is oriented.

Using the lemma of homogeneity and Corollary \ref{cor:1} there
exists a diffeotopy $\{d_t \colon t\in [0,1] \}$ of $M^n$ such
that $d_0 = id_{M^n}$, moreover

\begin{itemize}
\item $d_1(S(g)) = S(f)$ if $n$ is even,
\item $d_1(S_+(g)) = S_+(f)$ and $d_1(S_-(g)) = S_-(f)$ if $n$ is odd.
\end{itemize}

This implies that $S(f \circ d_1) = S(g)$, moreover $S_+(f \circ
d_1) = S_+(g)$ and $S_-(f \circ d_1) = S_-(g)$ if $n$ is odd.
Since $\{f \circ d_t \colon t \in [0,1]\}$ provides a regular
homotopy connecting $f$ and $f \circ d_1$ we might suppose that
$S(f) = S(g)$, moreover $S_+(f) = S_+(g)$ and $S_-(f) = S_-(g)$ if
$n$ is odd.

If $n$ is even then any two cross-caps are locally equivalent. For
$n$ odd any two cross-caps of the same sign are locally
equivalent. Thus there is a neighborhood $U$ of $S(f) = S(g)$ such
that $f$ is regularly homotopic to $g$ on $U$. This regular
homotopy can be extended to $M^n$ using Smale's lemma. Thus we may
suppose that $f|U = g|U$.

If $U$ is chosen small enough then $f|U = g|U$ is generic. We
perturb $f$ and $g$ slightly outside $U$, using a regular
homotopy, to obtain generic maps. Thus we can also suppose that
$f, g \in G(M^n, \Real^{2n-1})$.

Choose a smooth simplicial decomposition of $M^n$ so fine that for
any $n$-simplex $\Delta^n$ containing a point $p$ of $S(f) = S(g)$
we have $\Delta^n \subset U$ and $p \in \text{int}(\Delta^n)$.
(The proof of the existence of a smooth simplicial decomposition
can be found in \cite{Munkres}.) Thus $f| \Delta^n = g| \Delta^n$
for each $\Delta^n$ as above. For $p \in S(f)$ we also choose
$\Delta^n \ni p$ so small that $f|\Delta^n = g|\Delta^n$ has the
canonical form (\ref{eqn:1}) in a coordinate neighborhood
containing $\Delta^n$ and centered at $p$.

The fiber of the bundle $\text{MONO}(TM, T\Real^{2n-1})$ is
homeomorphic to the Stiefel manifold $V_{2n-1, n}$ of $n$-frames
in $\Real^{2n-1}$. It is well known that $\pi_i(V_{2n-1,n}) = 0$
for $i < n-1$.

Thus there exists an $M^n$-regular homotopy connecting the $M^n$-
immersion $$(f, df)|\text{sk}_{n-2}(M^n) \,\, \text{with} \,\, (g,
dg)|\text{sk}_{n-2}(M^n).$$ Using Smale's lemma we can extend this
$M^n$-regular homotopy to the whole manifold $M^n$ keeping $f|U$
fixed. So we might suppose that
$$ (f, df)|\text{sk}_{n-2}(M^n) = (g, dg)|\text{sk}_{n-2}(M^n). $$

\subsection{The first obstruction}

Our next task is to find an $M^n$-regular homotopy connecting $(f,
df)$ with $(g, dg)$ on $sk_{n-1}(M^n)$. The obstruction
$\Omega(df,dg)$ to finding a regular homotopy connecting $(f, df)$
with $(g, dg)$ on an $(n-1)$-simplex $\Delta^{n-1}$, fixing the
boundary $\partial \Delta^{n-1}$, lies in
\[
\pi_{n-1}(V_{2n-1,n}) \approx \left\{
\begin{array}{ll}
\mathbb{Z}_2 & \mbox{if $n$ is even,}\\
\mathbb{Z} & \mbox{if $n$ is odd.}
\end{array}
\right.
\]

We now prove that the sum of the obstructions on the boundary of
an $n$-simplex on which both $f$ and $g$ are immersions vanishes.

\begin{lem} \label{lem:6}
Suppose that the $C^{\infty}$ $\Real^n$-immersions $(f,f')$ and
$(g,g')$ of $\partial\Delta^n$ into $\Real^q$ $(q>n)$ are tangent
on $\text{sk}_{n-2}(\Delta^n)$. Denote the faces of $\Delta^n$ by
$$\Delta^{n-1}_0, \dots, \Delta^{n-1}_n,$$
oriented by the outward normal vectors of $\Delta^n$ and let
$$(f_i,f_i') = (f,f')|\Delta^{n-1}_i, \,\,\,\, (g_i,g_i') =
(g,g')|\Delta^{n-1}_i \,\, (i = 0, \dots, n).$$ Then
$$\sum_{i=0}^n \Omega(f_i', g_i') = \tau(f')-\tau(g'). $$
\end{lem}

\begin{proof}
We define a sequence of $C^{\infty}$ $\Real^n$-immersions
$$(h_i, h_i')\colon \partial \Delta^n \to \Real^{2n-1} \,\, (i = 0,
\dots, n+1)$$ such that $(h_0,h_0')=(f,f')$, $(h_{n+1},h_{n+1}') =
(g,g')$ and for $0 < i < n+1$ let
\[
(h_i, h_i')|\Delta^{n-1}_j = \left\{
\begin{array}{ll}
(g_j,g_j') & \mbox{if $j < i$,}\\
(f_j,f_j') & \mbox{if $i \le j \le n$.}
\end{array}
\right.
\]
Then $(h_i,h_i')$ is $C^{\infty}$ since $(f,f')$ and $(g,g')$ are
tangent on $sk_{n-2}(\Delta^n)$. Applying Theorem \ref{thm:3} to
$(h_i,h_i')$ and $(h_{i+1},h_{i+1}')$ we get that $$\Omega(f_i',
g_i') = \tau(h_i') - \tau(h_{i+1}').$$ Summing up these equations
we obtain the required result.
\end{proof}

\begin{cor} \label{cor:2}
Let $q>n$ and $f,g \in \imm(\Delta^n, \Real^q)$. Suppose that
$(f,f') = (f,df)|\partial \Delta^n$ and $(g,g') = (g, dg)|\partial
\Delta^n$ are tangent on $\text{sk}_{n-2}(\Delta^n)$. Using the
notations of Lemma \ref{lem:6} it holds that
$$\sum_{i=0}^n \Omega(f_i', g_i') = 0.$$
\end{cor}

\begin{proof}
Since both $(f,f')$ and $(g,g')$ are extendible we have that
$\tau(f') = 0$ and $\tau(g') = 0$. Now Lemma \ref{lem:6} gives the
statement of our corollary.
\end{proof}

Lemma \ref{lem:2} implies that there is a homeomorphism  $$\chi
\colon \mathcal{I}'(n,2n-1;n) \to \imm_1(S^{n-1}, \Real^{2n-1}).$$

\begin{note}
Let $\text{Emb}_1(S^k, \Real^l)$ denote the space of pairs $(h,
\nu)$, where $h \colon S^k \to \Real^l$ is an embedding and $\nu$
is a transversal $1$-field along $h$.
\end{note}

\begin{defn}
Two elements of $\imm_1(S^{n-1}, \Real^{2n-1})$ (or
$\mathcal{I}'(n,2n-1;n)$) are called \emph{regularly homotopic} if
they lie in the same path-component of the corresponding space. A
\emph{regular homotopy} is a path in $\imm_1(S^{n-1},
\Real^{2n-1})$ (or $\mathcal{I}'(n,2n-1;n)$).
\end{defn}

\begin{prop} \label{prop:5}
For every immersion $(f, \nu) \in \imm_1(S^{n-1},\Real^{2n-1})$
there exists an embedding $(h, \mu) \in
\text{Emb}_1(S^{n-1},\Real^{2n-1})$ such that $(f, \nu)$ is
regularly homotopic to $(h, \mu)$.
\end{prop}

\begin{proof}
The subspace of embeddings is dense open in $\imm(S^{n-1},
\Real^{2n-1})$. Thus we can perturb $f$ by a regular homotopy to
get an embedding $h$. If the perturbation is sufficiently small
then $\nu$ can be deformed simultaneously as a transversal field
along the regular homotopy.
\end{proof}

\begin{defn}
Suppose that $\Real^{2n-1}$ is oriented. Let
\[
\text{lk} \colon \text{Emb}_1(S^{n-1}, \Real^{2n-1}) \to \left\{
\begin{array}{ll}
\mathbb{Z} & \mbox{if $n$ is odd,}\\
\mathbb{Z}_2 & \mbox{if $n$ is even,}
\end{array}
\right.
\]
be defined as follows: For $(f, \nu) \in \text{Emb}_1(S^{n-1},
\Real^{2n-1})$ choose $\varepsilon > 0$ sufficiently small so that
$f + \varepsilon \nu$ is an embedding. Then let $\text{lk}(f,\nu)$
be accordingly the $\mathbb{Z}$ or $\mathbb{Z}_2$ linking number
$\text{lk}(f, f + \varepsilon \nu)$.
\end{defn}

\begin{prop} \label{prop:6}
$\text{lk}$ is a regular homotopy invariant, i.e., if the maps \\
$(f,\nu), (h, \mu) \in \text{Emb}_1(S^{n-1},\Real^{2n-1})$ are
regularly homotopic then $\text{lk}(f,\nu) = \text{lk}(h,\mu)$.
\end{prop}

\begin{proof}
This is exactly Lemma 9 in \cite{Whitney}.
\end{proof}

Now we can extend the definition of lk from embeddings to
immersions.

\begin{defn} \label{defn:2}
Suppose that $\Real^{2n-1}$ is oriented. Let
\[
\text{lk} \colon \imm_1(S^{n-1}, \Real^{2n-1}) \to \left\{
\begin{array}{ll}
\mathbb{Z} & \mbox{if $n$ is odd,}\\
\mathbb{Z}_2 & \mbox{if $n$ is even,}
\end{array}
\right.
\]
be defined as follows: For $(f, \nu) \in \imm_1(S^{n-1},
\Real^{2n-1})$ choose an embedding $(h,\mu)$ in
$\text{Emb}_1(S^{n-1},\Real^{2n-1})$ regularly homotopic to $(f,
\nu)$ (such an embedding exists because of Proposition
\ref{prop:5}). Then let $\text{lk}(f,\nu) = \text{lk}(h,\mu)$.
This definition does not depend on the choice of $(h, \mu)$
because of Proposition \ref{prop:6}.
\end{defn}

\begin{prop} \label{prop:1}
The function $\text{lk}$ defined above is a regular homotopy
invariant.
\end{prop}

\begin{proof}
This is clear from Proposition \ref{prop:6} and Definition
\ref{defn:2}.
\end{proof}

\begin{note} \label{note:2}
Let $I_n' = \pi_0(\mathcal{I}'(n,2n-1;n))$ and $J_n' =
\pi_0(\imm_1(S^{n-1}, \Real^{2n-1}))$. Proposition \ref{prop:1}
implies that $\text{lk}$ defines a function on $J_n'$, which we
denote by $\text{lk}_*$.
\end{note}

\begin{lem} \label{lem:7}
Let $f \in L(D^n, \Real^{2n-1})$ and $$(f,f') = (f,df)|S^{n-1} \in
\mathcal{I}'(n,2n-1;n).$$ Then $$\text{lk}(\chi(f,f')) = \#S(f).$$
(For the definition of $\chi$ see Lemma \ref{lem:2}.)
\end{lem}

\begin{proof}
Perturb $f$ using a regular homotopy to get a generic map $h$
(thus $h|S^{n-1}$ is an embedding). Let $(h,h') = (h,
dh)|S^{n-1}$. Proposition \ref{prop:1} implies that
$\text{lk}(\chi(f,f')) = \text{lk}(\chi(h,h'))$. Moreover, the
number of singular points is also a regular homotopy invariant,
thus $\#S(f) = \#S(h)$. By definition $\text{lk}(\chi(h,h')) =
\mathcal{L}_h(D^n)$. On the other hand Lemma \ref{lem:5} implies
that $\mathcal{L}_h(D^n) = \#S(h)$.
\end{proof}

\begin{cor} \label{cor:4}
Let $(f, \nu) \in \imm_1(S^{n-1}, \Real^{2n-1})$. Then
$$ \text{lk}(f,-\nu) = -\text{lk}(f, \nu).$$
\end{cor}

\begin{proof}
Choose two generic maps $f_0, f_1 \in G(D^n, \Real^{2n-1})$ such
that
$$\chi((f_i,df_i)|S^{n-1}) = (f, (-1)^i \nu)$$ for $i=0,1$. Then
$f_0$ and $f_1$ fit together to define a map $h \in G(S^n,
\Real^{2n-1})$. Lemma \ref{lem:5} implies that $\#S(f_0)+ \#S(f_1)
= \#S(h) = 0$. Moreover, Lemma \ref{lem:7} implies that
$\text{lk}(f, \nu) = \#S(f_0)$ and $\text{lk}(f, -\nu) =
\#S(f_1)$. Putting together these results we get that
$\text{lk}(f, \nu) + \text{lk}(f, -\nu) = 0$, as required.
\end{proof}

\begin{prop} \label{prop:2}
$\tau$ (see Definition \ref{defn:3}) is a regular homotopy
invariant. Moreover, it induces a bijection $\tau_* \colon I_n'
\to \pi_{n-1}(V_{2n-1,n})$.
\end{prop}

\begin{proof}
A regular homotopy of a map $(f,f') \in \mathcal{I}'(n,2n-1;n)$
induces a homotopy of the map  $S^{n-1} \to V_{q,n}$ defined by
$$x \mapsto f'\{e_1(x), \dots, e_n(x)\}.$$ Thus the homotopy class
$\tau(f')$ is constant on the regular homotopy class of $(f,f')$.

Hirsch's theorem (Theorem \ref{thm:4}) implies that $$\tau_*
\colon \pi_0(\imm_1(S^{n-1}, \Real^{2n-1})) \to
\pi_0(\text{MONO}(TS^{n-1} \oplus \varepsilon^1, \Real^{2n-1}))
\approx \pi_{n-1}(V_{2n-1,n})
$$ is a bijection (the bundle $TS^{n-1} \oplus \varepsilon^1$ is trivial).
\end{proof}

\begin{defn}
Define the \emph{connected sum} operation $\#$ on elements $x$ and
$y$ of $\imm_1(S^{n-1}, \Real^{2n-1})$ or $\mathcal{I}'(n,2n-1;n)$
by joining them with a thin tube. This operation is also well
defined on regular homotopy classes of maps: Suppose that $x
\sim_H x_1$ and $y \sim_K y_1$. Smale's lemma implies that the
regular homotopies $H$ and $K$ may be kept fixed on small disks.
Now a regular homotopy $L$ connecting $x \# x_1$ with $y \# y_1$
is defined by taking $H$ and $K$, and joining them with a tube
attached to the disks kept fixed.
\end{defn}

\begin{prop} \label{prop:3}
The sets $I_n'$ and $J_n'$ (introduced in Notation \ref{note:2})
endowed with the connected sum operation form abelian groups.
Moreover, $\chi_* \colon I_n' \to J_n'$ and $\tau_*$ are group
isomorphisms and $\text{lk}_*$ is a group homomorphism.
\end{prop}

\begin{proof}
It is clear that $(I_n', \#)$ and $(J_n',\#)$ are abelian
semigroups. Let $i \colon D^n \to \Real^{2n-1}$ denote the
standard embedding. Then $(i,i')=(i,di)|S^{n-1}$ (respectively \\
$\chi((i,di)|S^{n-1}) = (i, r)$, where $r$ is the outward normal
field of $S^{n-1}$ in $\Real^n$) represents the identity of $I_n'$
(respectively $J_n'$).

If $x,y \in I_n'$ then choose a representative $(f,f')$ for $x$
and $(g,g')$ for $y$. Then $\tau(f'\# g')$ is the connected sum of
the spheroids $\tau(f')$ and $\tau(g')$, which implies that
$\tau_*(x \# y) = \tau_*(x) + \tau_*(y)$. Since $(i,i')$ is
extendible, $\tau_*(i') = 0$. Proposition \ref{prop:2} states that
$\tau_*$ is a bijection. Thus $\tau_*$ is a semigroup isomorphism
that takes the identity of $I_n'$ to the identity of
$\pi_{n-1}(V_{2n-1,n}) \approx \mathbb{Z}$ or $\mathbb{Z}_2$. This
implies that $I_n'$ is a group isomorphic to $\mathbb{Z}$ or
$\mathbb{Z}_2$ (depending on the parity of $n$) and $\tau_*$ is a
group isomorphism.

Since $\chi$ is a homeomorphism, $\chi_* \colon I_n' \to J_n'$ is
a bijection. Moreover, from the definition of $\chi$ it is easy to
see that $\chi_*$ is a semigroup homomorphism. We saw above that
$\chi_*$ takes the identity $(i,i')$ of $I_n'$ to the identity
$(i,r)$ of $J_n'$. Thus $J_n'$ is also a group isomorphic to
$\mathbb{Z}$ or $\mathbb{Z}_2$ and $\chi_*$ is a group
isomorphism.

We finally show that $\text{lk}_*$ is a group homomorphism from
$J_n'$ to $\mathbb{Z}$ or $\mathbb{Z}_2$. Suppose that $x,y \in
J_n'$ and choose representatives $(f_1, \nu_1)$ and $(f_2, \nu_2)$
for $x$, respectively $y$, such that $f_1$ and $f_2$ can be
separated by a hyperplane in $\Real^{2n-1}$. Also choose $h_1,h_2
\in G(D^n, \Real^{2n-1})$ so that $\chi((h_i,dh_i)|S^{n-1}) =
(f_i,\nu_i)$ for $i=1,2$. Denote the boundary sum $h_1 \natural
h_2$ with $h$. Then $\chi((h,dh)|S^{n-1}) = (f_1, \nu_1) \# (f_2,
\nu_2)$ and $\#S(h) = \#S(h_1) + \#S(h_2)$. Using Lemma
\ref{lem:7} we get that $\text{lk}((f_1,\nu_1)\#(f_2,\nu_2)) =
\text{lk}(f_1, \nu_1) + \text{lk}(f_2, \nu_2)$, i.e., $\text{lk}(x
\# y ) = \text{lk}(x)+ \text{lk}(y)$.
\end{proof}

The following lemma is the key to connecting the works of Hirsch
and Whitney.

\begin{lem} \label{lem:8}
There exists an isomorphism
\[
\alpha \colon \pi_{n-1}(V_{2n-1,n}) \to \left\{
\begin{array}{ll}
\mathbb{Z} & \mbox{if $n$ is odd,}\\
\mathbb{Z}_2 & \mbox{if $n$ is even,}
\end{array}
\right.
\]
such that
$$\text{lk} \circ \chi = \alpha \circ \tau. $$ I.e., the following
diagram is commutative:
\[
\begin{CD}
\mathcal{I}'(n,2n-1;n) @>\chi>>
\imm_1(S^{n-1}, \Real^{2n-1}) \\
@V\tau VV @V\text{lk}VV \\
\pi_{n-1}(V_{2n-1,n}) @>\alpha >> \mathbb{Z} \,\,\,\text{or}\,\,\,
\mathbb{Z}_2.
\end{CD}
\]
\end{lem}

\begin{rem}
If $n$ is even then $\alpha$ might be omitted from the formula
since $\mathbb{Z}_2$ has only one automorphism.
\end{rem}

\begin{proof}
We saw that $\tau$ and $\text{lk}$ are regular homotopy
invariants, i.e., they are constant on the path-components of
their domains. They define maps $$\tau_* \colon  I'_n =
\pi_0(\mathcal{I}'(n,2n-1;n)) \to \pi_{n-1}(V _{2n-1,n})$$ and
\[
\text{lk}_* \colon J'_n= \pi_0(\imm_1(S^{n-1}, \Real^{2n-1})) \to
\left\{
\begin{array}{ll}
\mathbb{Z} & \mbox{if $n$ is odd,}\\
\mathbb{Z}_2 & \mbox{if $n$ is even.}
\end{array}
\right.
\]
Proposition \ref{prop:3} implies that $I'_n$ and $J'_n$ are
abelian groups with the connected sum operation, moreover,
$\chi_*, \tau_*$ are isomorphisms  and $\text{lk}_*$ is a
homomorphism. We also know that $I_n'$ and $J_n'$ are isomorphic
to $\mathbb{Z}$ or $\mathbb{Z}_2$, according to the parity of $n$.

$\text{lk}_*$ is an epimorphism because of Lemma \ref{lem:4}:
$$\text{lk}(\chi(w,w'))= \pm 1.$$ So $\text{lk}_*$ is a
$\mathbb{Z} \to \mathbb{Z}$ or $\mathbb{Z}_2 \to \mathbb{Z}_2$
epimorphism, and thus it is also an isomorphism. Now define
$\alpha$ to be  $\text{lk}_* \circ \chi_* \circ (\tau_*)^{-1}$.
\end{proof}

From Theorem \ref{thm:3} and Lemma \ref{lem:8} we get the
following

\begin{prop}
Suppose that $(f,\nu), (g,\mu) \in \imm_1(\partial\Delta^n,
\Real^q)$ agree on $(\partial\Delta^n) \setminus
\text{int}(\Delta^{n-1})$, where $\Delta^{n-1}$ is an $(n-1)$-face
of the $n$-simplex $\Delta^n$. Then for $(f,f') =
(\chi^{-1}(f,\nu))|\Delta^{n-1}$ and $(g,g') =
(\chi^{-1}(g,\mu))|\Delta^{n-1}$
$$\alpha(\Omega(f', g')) = \text{lk}(f,\nu)-\text{lk}(g,\mu).$$
\end{prop}

\begin{rem}
$\Delta^{n-1}$ is co-oriented by the outward normal vector of
$\Delta^n$. Then $$(\chi^{-1}(f,\nu))|\Delta^{n-1} =
\zeta^{-1}((f, \nu)|\Delta^{n-1}),$$ where $\zeta$ is the
homeomorphism of Lemma \ref{lem:1}. Thus
$\text{lk}(f,\nu)-\text{lk}(g,\mu)$ depends only on
$(f,\nu)|\Delta^{n-1}$ and $(g,\mu)|\Delta^{n-1}$.
\end{rem}

\begin{defn}
For $(f,f'),(g,g') \in \mathcal{I}(n-1,2n-1;n)$ that are tangent
on $\partial\Delta^n$ let
$$o(f',g') = \alpha(\Omega(f',g')).$$
\end{defn}

\begin{cor} \label{cor:5}
Let $\Delta^{n-1}$ be a face of the $n$-simplex $\Delta^n$ and let
$(f,f')$ and $(g,g')$ be $C^{\infty}$ $\Real^n$-immersions of
$\Delta^{n-1}$ into $\Real^{2n-1}$ tangent on
$\partial\Delta^{n-1}$. Choose extensions of $(f,f')$ and $(g,g')$
to $\partial\Delta^n$ that are tangent on $\partial\Delta^n
\setminus \Delta^{n-1}$. Then $$o(f',g') = \text{lk} \circ
\chi(f,f') - \text{lk} \circ \chi(g,g').$$
\end{cor}

\begin{lem} \label{lem:10}
Let $(f,f')$ and $(g,g')$ be as in Corollary \ref{cor:5} and
$(w,w')$ be the boundary of a Whitney-umbrella of sign
$\varepsilon$ (as in Notation \ref{note:1}). Then $$o(f'\# w', g')
= o(f',g') + \varepsilon.$$
\end{lem}

\begin{proof}
Using extensions as in Corollary \ref{cor:5} we get that $$o(f'\#
w', g') = \text{lk} \circ \chi((f,f')\#(w,w')) - \text{lk} \circ
\chi(g,g').$$ Using the fact that $\text{lk}$ and $\chi$ are
additive we get that this equals $$\text{lk} \circ \chi(f,f')+
\text{lk} \circ \chi(w,w') - \text{lk} \circ \chi(g,g').$$ The
definition of the sign of an umbrella implies that $\text{lk}
\circ \chi(w,w') = \varepsilon $.
\end{proof}

\begin{prop} \label{prop:4}
Let $(f,f')$ and $(g,g')$ be as in Corollary \ref{cor:5}. Then
$(f,f')$ and $(g,g')$ are regularly homotopic (rel
$\partial\Delta^{n-1}$) iff $o(f',g')=0$.
\end{prop}

\begin{proof}
This is a trivial consequence of Lemma \ref{lem:9}.
\end{proof}

\begin{rem}
Lemma \ref{lem:10} and Proposition \ref{prop:4} show us how
Whitney-umbrellas destroy obstructions to moving $f$ to $g$ on an
$(n-1)$-simplex (rel boundary). We have to take the connected sum
of $f$ with several copies of the boundary of an umbrella. This
will be done by a diffeotopy of $M^n$ that pushes the
$(n-1)$-simplex through a singular point.
\end{rem}

\begin{lem} \label{lem:11}
Using the notations of Lemma \ref{lem:6}
$$\sum_{i=0}^n
o(f_i', g_i') = \text{lk}\circ\chi(f,f')-\text{lk}\circ\chi(g,g').
$$  In particular, if $(f,f')$ and $(g,g')$ are both extendible
then $$\sum_{i=0}^n o(f_i',g_i') = 0.$$
\end{lem}

\begin{proof}
This is a trivial consequence of Lemma \ref{lem:6}.
\end{proof}

\begin{rem}
Lemma \ref{lem:11} shows us that the obstruction $o$ has the
cocycle property, i.e., the sum of the obstructions on the
boundary of an $n$-simplex, on which both $f$ and $g$ are
immersions, vanishes.
\end{rem}

\subsection{Pushing (n-1)-simplices through Whitney-umbrellas}

Now we are going to use the apparatus developed above for the
proof of Theorem \ref{thm:1}. Since until now we only worked in a
standard simplex in Euclidean space, we have to globalize our
results to the triangulated manifold $M^n$. We do \emph{not}
require $M^n$ to be oriented.

\begin{defn}
Suppose that $\Delta^n$ is an $n$-simplex of $M^n$. Then we can
define the map $\chi = \chi_{\Delta^n}$ taking $M^n$-immersions of
$\partial\Delta^n$ to $\imm_1(\partial\Delta^n, \Real^{2n-1})$
using the outward normal vectors $r(x)$ of $\partial\Delta^n$
(corners rounded off) as follows: Let $(f,f')$ be an
$M^n$-immersion of $\partial\Delta^n$ into $\Real^{2n-1}$. Then
let $\chi(f,f') = (f, \nu)$, where $\nu(x) = f'(r(x))$.
\end{defn}

\begin{defn}
Suppose that $\Delta^{n-1}$ is an $(n-1)$-simplex of $M^n$ that is
co-oriented. Then we can define the map $\zeta =
\zeta_{\Delta^{n-1}}$ taking $M^n$-immersions of $\Delta^{n-1}$ to
$\imm_1(\Delta^{n-1}, \Real^{2n-1})$ as follows: Let $r(x)$ be a
positive normal field along $\Delta^{n-1}$. If $(f,f')$ is an
$M^n$-immersion of $\Delta^{n-1}$ into $\Real^{2n-1}$ then let
$\zeta(f,f') = (f, \nu)$, where $\nu(x) = f'(r(x))$.
\end{defn}

Corollary \ref{cor:5} implies that the following definition makes
sense:

\begin{defn}
Suppose that $\Delta^{n-1}$ is an $(n-1)$-simplex of $M^n$ that is
co-oriented and let $(f,f')$ and $(g,g')$ be $M^n$-immersions of
$\Delta^{n-1}$ into $\Real^{2n-1}$ tangent on
$\partial\Delta^{n-1}$. Then the obstruction $o(f',g')$ to finding
a regular homotopy of $(f,f')$ to $(g,g')$ (rel
$\partial\Delta^{n-1}$) is given as follows: Choose an $n$-simplex
$\Delta^n$ of $M^n$ such that $\Delta^{n-1} \subset
\partial\Delta^n$. Extend $(f,f')$ and $(g,g')$ to
$\partial\Delta^{n}$ to be tangent on $\partial\Delta^n \setminus
\Delta^{n-1}$. Also choose a normal field $s(x)$ along
$\partial\Delta^n$ that agrees with the co-orientation of
$\Delta^{n-1}$. Now define $o(f',g')$ to be $\text{lk}(f,f' \circ
s) - \text{lk}(g,g' \circ s)$.
\end{defn}

\begin{rem} \label{rem:2}
If we change the co-orientation of $\Delta^{n-1}$ then $o(f',g')$
changes sign: To see this we choose the same extensions of
$(f,f')$ and $(g,g')$ to the boundary $\partial \Delta^n$ of the
same $n$-simplex $\Delta^n$. Then we must take $-s$ instead of $s$
because of the changed co-orientation of $\Delta^n$. Now using
Corollary \ref{cor:4} we get that
$$\text{lk}(f,f' \circ (-s)) - \text{lk}(g,g' \circ (-s)) =
-\text{lk}(f,f' \circ s)) + \text{lk}(g,g' \circ s) = -o(f',g').$$
\end{rem}

Let $f$ and $g$ be the two locally generic maps of Theorem
\ref{thm:1}. We have supposed that $(f,df)$ and $(g,dg)$ are
tangent on $\text{sk}_{n-2}M^n$. If $n$ is odd choose a
co-orientation $O_{\Delta^{n-1}}$ for every $(n-1)$-simplex
$\Delta^{n-1}$ of $M^n$ so that $o(df,dg) \le 0$. This can be done
because of Remark \ref{rem:2}. Let $k = \abs{S(f)}=\abs{S(g)}$ and
denote by $\Delta^n_1, \dots, \Delta^n_k$ the $n$-simplices of
$M^n$ that contain any point of $S(f)=S(g)$.

We are now going to construct an oriented curve $\gamma$ on $M^n$
that intersects each simplex $\Delta^{n-1}$ of $M^n$ in the
positive direction (according to $O_{\Delta^{n-1}}$) in exactly
$\abs{o(df,dg)}$ points. Notice that if $\Delta^{n-1} \subset
\partial\Delta^n_i$ for some $1 \le i \le k$ then on $\Delta^{n-1}$ we
have $o(df,dg)=0$ since $f|\Delta^n_i = g|\Delta^n_i$. Thus
$\gamma$ avoids the simplices $\Delta^n_i$. This $\gamma$ may be
thought of as the dual of the obstruction to deform $f$ to $g$ on
the $(n-1)$-skeleton of $M^n$.

To obtain $\gamma$ choose a set of $\abs{o(df,dg)}$ points
$P_{\Delta^{n-1}}$ on each $(n-1)$-simplex of $M^n$. Now take an
$n$-simplex $\Delta^n$ and denote its faces by $\Delta^{n-1}_0,
\dots, \Delta^{n-1}_n$, moreover let $(f_i,f_i') =
(f,df)|\Delta^{n-1}_i$, $(g_i,g_i') = (g,dg)|\Delta^{n-1}_i$, $O_i
= O_{\Delta^{n-1}_i}$ and $P_i = P_{\Delta^{n-1}_i}$ for $0 \le i
\le n$. Lemma \ref{lem:11} implies that if each $\Delta^{n-1}_i$
is co-oriented by the outward normal vector of $\Delta^n$ (denote
this co-orientation by $U_i$) then
\begin{equation} \label{eqn:2}
\sum_{i=0}^n o(f_i',g_i') = 0.
\end{equation}
If we consider $\Delta^{n-1}_i$ with the co-orientation $U_i$ then
$o(f_i',g_i') < 0$ implies that $U_i = O_i$ and $o(f_i',g_i')
> 0$ implies that $U_i = -O_i$. Give a minus sign to each
point of $P_i$ if $U_i = O_i$ and a plus sign otherwise. Then
equation $(\ref{eqn:2})$ is equivalent to the statement that the
sum of the signs in $\bigcup_{i=0}^n P_i$ equals $0$. Now let
$\gamma \cap \Delta^n$ be given as follows: Make a bijection
between the $+$ and $-$ points of $\bigcup_{i=0}^n P_i$ and
connect each pair of points with an embedded curve segment
oriented from $+$ to $-$. We do this so that these curve segments
are pairwise disjoint. This is possible since $n \ge 3$.

Doing this for each $n$-simplex we obtain an oriented embedded
curve $\gamma$ with the required intersection property. Now we
make $\gamma$ connected by taking the connected sum of its
components: Let for example $\gamma_1$ and $\gamma_2$ be two
components of $\gamma$. Then choose two points $p_1 \in \gamma_1$
and $p_2 \in \gamma_2$. Join $p_1$ and $p_2$ with an embedded
curve $\eta$ that avoids $\Delta^n_i$ for $1 \le i \le k$ and also
each $(n-2)$-simplex and intersects each $(n-1)$-simplex
transversally. Then take two parallel curves $\eta_1$ and $\eta_2$
close to $\eta$ and orient them according to the orientations of
$\gamma_1$ and $\gamma_2$. If $\eta$ intersects an $(n-1)$-simplex
$\Delta^{n-1}$ at a point $x$ then $\eta_1$ and $\eta_2$ will
intersect $\Delta^{n-1}$ in different directions near $x$.

Similarly, we might modify $\gamma$ so that it goes through
exactly one point $p$ of $S(f)$, and let the Whitney-umbrella of
$f$ at $p$ be positive if $n$ is odd. To do this, choose a small
embedded curve containing $p$ and join it to $\gamma$ as above. We
will also call this modified curve $\gamma$. Let us suppose that
the $n$-simplex containing $p$ is $\Delta^n_1$.

Thus it will hold true that the connected embedded curve $\gamma$
intersects \emph{algebraically}  each $(n-1)$-simplex
$\Delta^{n-1}$ in $\abs{o(df,dg)}$ points if we consider
$\Delta^{n-1}$ with the co-orientation $O_{\Delta^{n-1}}$.
Moreover, $\gamma$ contains the cross-cap point $p \in \Delta^n_1$
and is disjoint from $\Delta^n_i$ for $1 < i \le k$.

Let $\nu_{\gamma}$ be a thin tubular neighborhood of $\gamma$ in
$M^n$. We are going to construct a diffeotopy $\{\Psi_t \colon 0
\le t \le 1 \}$ of $\nu_{\gamma}$ such that $\Psi_0 =
id_{\nu_{\gamma}}$ and for every $t \in [0,1]$ the diffeomorphism
$\Psi_t$ is the identity in a neighborhood of $\partial
\nu_{\gamma}$. Thus we can extend $\Psi_t$ to the whole manifold
$M^n$ to be the identity outside $\nu_{\gamma}$.

\begin{figure}[t]
\includegraphics{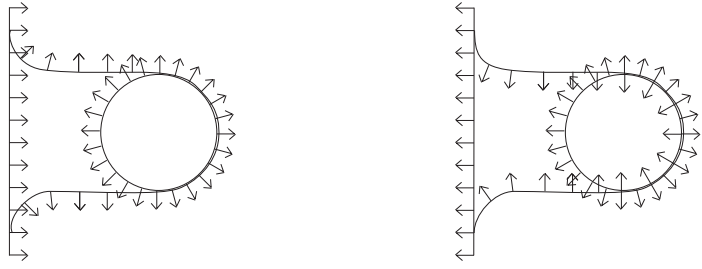}
\caption{Co-orientations and the connected sum} \label{fig:2}
\end{figure}

The diffeomorphism $\Psi_t$ is constructed as follows: Denote by
$T$ the identity of $\Real^{n-1}$ if $\gamma$ is orientation
preserving, and let $T$ be a reflection in a hyperplane of
$\Real^{n-1}$ if $\gamma$ is orientation reversing. Then
$\nu_{\gamma}$ is diffeomorphic to the factor space $$\Gamma =
\Real \times D^{n-1}_2 \left/_{(x,y) \sim (x+1, T(y))}\right.,$$
where $D^{n-1}_2 \subset \Real^{n-1}$ is the disc of radius $2$.
Denote by $p \colon \Real \times D^{n-1}_2 \to \Gamma$ the
projection. Let $\lambda \colon \Real \to \Real$ be a $C^{\infty}$
function such that $\lambda(x) = 0$ for $x < \varepsilon$ and
$\lambda(x) = 1$ for $x > 1- \varepsilon$ (where $\varepsilon < 1$
is a small positive constant). First we define a diffeotopy
$\{\phi_t \colon 0 \le t \le 1\}$ of $\Real \times D^{n-1}_2$ with
the formula
\[
\phi_t(x,y) = \left\{
\begin{array}{ll}
(x+t, y) & \mbox{if $\norm{y} \le 1$,}\\
(x + \lambda(s)t, y) & \mbox{if $\norm{y} = 2-s$.}
\end{array}
\right.
\]
$(x,y) \sim (x+1, T(y))$ implies that $(x+c,y) \sim (x+c+1, T(y))$
for any number $c$. Thus the diffeomorphism $\phi_t$ factors
through the projection $p$ to a diffeomorphism  $\Psi_t$ of
$\Gamma \approx \nu_{\gamma}$. I.e., the following diagram is
commutative:
\[
\begin{CD}
\Real\times D^{n-1}_2 @>\phi_t>>
\Real\times D_2^{n-1} \\
@VpVV @VpVV \\
\nu_{\gamma} @>\Psi_t >> \nu_{\gamma}.
\end{CD}
\]
Since $\Psi_1$ is the identity on $p(\Real \times D^{n-1}_1)$, the
diffeomorphism $\Psi_1$ is the identity on a thinner tubular
neighborhood of $\gamma$. Thus $\Psi_1(p) = p$. Moreover, let $F$
denote a fiber of $\nu_{\gamma}$ (diffeomorphic to $D^{n-1}$) with
a normal framing $v$ (in $M^n$) in the direction determined by
$\gamma$; and $S$ a small sphere around $p$ contained in
$\nu_{\gamma}$ with the outward normal framing $r$. Then the
framed submanifold $(\Psi_1(F),d\Psi_1(v))$ is equal to the
connected sum $(F,v) \# (S,r)$. On the other hand, the framed
submanifold $(\Psi_1(F), d\Psi_1(-v)) = (F,-v) \# (S,-r)$ (see
Figure \ref{fig:2}).

Now look at the regular homotopy $f \circ \Psi_t$ connecting $f$
with $h = f \circ \Psi_1$. Then $S(h) = S(f)$, moreover
$S_{\pm}(h) = S_{\pm}(f)$ if $n$ is odd. We also get that
$$(h,dh)|F = ((f,df)|F) \# ((f,df)|S).$$ Notice that $(f,df)|S$ is
the $M^n$-immersion $(w,w')$ of Notation \ref{note:1} and if we
co-orient $S$ by $r$ then $\text{lk} \circ \chi(w,w') = 1$ since
the sign of $p$ is $1$. Thus if we co-orient $F$ by $v$, as above,
then Lemma \ref{lem:10} gives that $o(dh|F, df|F) = 1$. On the
other hand, if we co-orient $F$ by $-v$ then $o(dh|F, df|F) = -1$
(see Figure \ref{fig:2}). So if we take an $(n-1)$-simplex
$\Delta^{n-1}$ of $M^n$ then $o(dh|\Delta^{n-1}, df|\Delta^{n-1})$
is the algebraic intersection of $\Delta^{n-1}$ and $\gamma$. From
the construction of $\gamma$ we know that
$$\Delta^{n-1}\cap \gamma = \abs{o(df|\Delta^{n-1},
dg|\Delta^{n-1})}.$$ Thus on $\Delta^{n-1}$ $$o(dh,dg) = o(dh,df)
+ o(df,dg) = \abs{o(df,dg)} + o(df,dg) = 0,$$ since the
co-orientation $O_{\Delta^{n-1}}$ was chosen so that $o(df,dg) \le
0$.

Consequently, there is no obstruction to finding an $M^n$-regular
homotopy between $h$ and $g$ on the $(n-1)$-skeleton of $M^n$.
Since $f$ is regularly homotopic to $h$, we suppose from now on
that $f = h$. Therefore we can suppose that the maps $f$ and $g$
coincide on (a neighborhood of) $\text{sk}_{n-1}(M^n)$.

\subsection{The second obstruction}

Let us recall that $\Delta^n_1, \dots, \Delta^n_k$ denote those
$n$-simplices of $M^n$ which contain a singular point of $f$,
where $k = 2l$ is even. (Each $\text{int}(\Delta^n_i)$ contains a
single Whitney umbrella point. If $n$ is odd then $\Delta^n_1,
\dots, \Delta^n_l$ contain a positive Whitney umbrella,
$\Delta^n_{l+1}, \dots, \Delta^n_{2l}$ a negative one for both $f$
and $g$. Moreover, $f|\Delta^n_i = g|\Delta^n_i$ for $1 \le i \le
k$.) Let
$$ U = \bigcup_{i=1}^k \text{int}(\Delta^n_i).$$ Take the manifold
with boundary $ N^n = M^n \setminus U $ and let $\hat{f} = f|N^n$
and $\hat{g} = g|N^n$. We know from the previous section that
$\hat{f}$ and $\hat{g}$ are regularly homotopic on
$\text{sk}_{n-1}(N^n)$. The secondary obstruction to finding a
regular homotopy between $\hat{f}$ and $\hat{g}$, not fixing the
boundary, lies in the cohomology group $H^n(N^n;
\pi_n(V_{2n-1,n}))$ with twisted coefficients. However, this group
is $0$ since an $n$-manifold without closed components is homotopy
equivalent to an $(n-1)$-dimensional simplicial complex.
% (Alternatively, this can be seen from the following short exact
% sequence:
% $$0 \longrightarrow H^n(N^n) \otimes G \longrightarrow H^n(N^n, G)
% \longrightarrow \text{Tor}(H^{n+1}(N^n), G) \longrightarrow 0,$$
% where $G = \pi_n(V_{2n-1,n})$. Since $\partial N^n \ne \emptyset$,
% we get that $H^n(N^n) = 0$ and thus $H^n(N^n; G) = 0$.)

So there exists a regular homotopy connecting $\hat{f}$ with
$\hat{g}$, but not necessarily fixing $\partial N^n$. Using
Smale's lemma we can extend this regular homotopy to the whole
manifold $M^n$, fixing a small neighborhood $V$ ($\subset U$) of
$S(f) = S(g)$. (During this regular homotopy on the "collar"
$\overline{U} \setminus V$ the map $f$ might become twisted, i.e.,
it may differ from $g$.)

Hence we might suppose that $f|N^n = g|N^n$ and $f|V = g|V$, but
it might happen that $f|(U \setminus V) \neq g|(U \setminus V)$.
If $n$ is odd then $\abs{S_+(f)} = \abs{S_+(g)}=l$, let us suppose
that
$$ S_+(f) = S_+(g) \subset \bigcup_{1 \le i \le l} \Delta^n_i, $$
while $$ S_-(f) = S_-(g) \subset \bigcup_{l+1 \le i \le k}
\Delta^n_i. $$

Since up to this point we have not used the assumption $n > 3$, we
can formulate and prove the following result that is weaker than
Theorem \ref{thm:1} and will not be used in its proof, but will be
needed in dealing with the exceptional case $n=3$.

\begin{thm}
Suppose that $n > 2$ and let $M^n$ be a closed $n$-manifold.
Moreover, let $f,g \in L'(M^n, \Real^{2n-1})$ be two locally
generic maps that satisfy $\abs{S(f)} = \abs{S(g)}$. Then there
exists a map $f' \in L'(M^n, \Real^{2n-1})$ and a subset $D
\subset M^n$ diffeomorphic to the disk $D^n$ such that $f \sim
f'$, the maps $f'|D$ and $g|D$ are immersions, and $f'|(M^n
\setminus D) = g|(M^n \setminus D)$.
\end{thm}

\begin{proof}
Let $U$ and $V$ be as above (where $S(g) \subset V \subset U$).
I.e., $U$ is the union of the interiors of the simplices
$\Delta^n_i$ for $1 \le i \le k$ (where $k = \abs{S(f)} =
\abs{S(g)}$) and $f$ is regularly homotopic to a map $f'$ such
that $f'|(M^n \setminus U) = g|(M^n \setminus U)$ and $f'|V =
g|V$. We can also suppose that $V_i = \Delta^n_i \cap V$ is
diffeomorphic to an $n$-disc for $1 \le i \le k$. Let us take a
small neighborhood $U_i$ of $\Delta^n_i$ that is diffeomorphic to
$D^n$. From now on we will denote by $U$ the union $\cup_{i=1}^k
U_i$. For each $i$ choose a $1$-simplex $e_i \subset U_i \setminus
V_i$ connecting a point of $\partial U_i$ with a point of
$\partial V_i$ and let $E_i$ be a thin tubular neighborhood of
$e_i$. Then $D_i = U_i \setminus (V_i \cup E_i)$ is diffeomorphic
to $D^n$. Since $n > 2$ there is no obstruction to constructing a
regular homotopy between $f'|E_i$ and $g|E_i$ fixing the ends of
the tube $E_i$. Thus (using Smale's lemma) we can suppose that
$f'|E_i = g|E_i$ for $1 \le i \le k$. Now we connect $D_i$ and
$D_{i+1}$ for $1 \le i \le k-1$ with a tube $T_i \subset M^n
\setminus U$ diffeomorphic to $I \times D^{n-1}$. (Thus each $T_i$
is an $n$-dimensional $1$-handle attached to $D_i$ and $D_{i+1}$.)
Let $$D = \left(\bigcup_{i=1}^k D_i \right) \cup
\left(\bigcup_{j=1}^{k-1} T_j \right).$$ Then $D$ is diffeomorphic
to $D^n$ and clearly $f'|D$ and $g|D$ are immersions, finally
$f'|(M^n \setminus D) = g|(M^n \setminus D)$.
\end{proof}

Let $F_1, \dots, F_l$ be pairwise disjoint subsets of $M^n$ such
that $F_i$ is diffeomorphic to the closed $n$-disc $D^n$, moreover
$F_i$ contains $\Delta^n_i$ and $\Delta^n_{l+i}$ in its interior
for $1 \le i \le l$. We will show that for $1 \le i \le l$ the
locally generic maps $f_i = f|F_i$ and $g_i = g|F_i$ are regularly
homotopic keeping a neighborhood of $\partial F_i$ fixed. This
would finish the proof of Theorem \ref{thm:1}. From now on let $1
\le i \le l$ be fixed.

It is apparent from the choice of $F_i$ that $S(f_i) = S(g_i)$ and
$\abs{S(f_i)} = \abs{S(g_i)}=2$. Moreover, if $n$ is odd we have
that $S_+(f_i) = S_+(g_i)$ and $S_-(f_i) = S_-(g_i)$ are both
$1$-element sets. Thus $\# S(f_i) = 0$ and $\# S(g_i) = 0$. Lemma
\ref{lem:5} implies that $\mathcal{L}_{f_i}(F_i) = 0$ and
$\mathcal{L}_{g_i}(F_i) = 0$. Using Lemma \ref{lem:8} we get that
$$0 = \mathcal{L}_{f_i}(F_i) = \text{lk} \circ \chi ((f_i,df_i)|\partial F_i) =
\alpha \circ \tau((f_i,df_i)|\partial F_i).$$ Since $\alpha$ is an
isomorphism, $\tau((f_i,df_i)|\partial F_i) = 0$ and from Lemma
\ref{lem:12} we can see that the map $(f_i,df_i)|\partial F_i$ is
extendible. Similarly, $(g_i,dg_i)|\partial F_i$ is also
extendible. (This also follows from Theorem \ref{thm:8}.)

Let $S^n_+$ (respectively $S^n_-$) denote the northern
(respectively southern) hemisphere of $S^n$ and we identify
$S^n_+$ with $F_i$. Then we get from the previous paragraph that
there exist two locally generic maps $f_i', g_i' \in L(S^n,
\Real^{2n-1})$ such that $f_i'|S^n_+ = f_i$, $g_i'|S^n_+ = g_i$,
moreover $f_i'|S^n_- = g_i'|S^n_-$ is an immersion. Thus it is
sufficient to prove the following theorem:

\begin{thm} \label{thm:6}
Let $f, g \in L(S^n, \Real^{2n-1})$ be locally generic maps and
suppose that $f|S^n_- = g|S^n_-$ are immersions, $\abs{S(f)} =
\abs{S(g)} = 2$ and $S(f) = S(g)$, moreover $S_+(f) = S_+(g)$ if
$n$ is odd. Then there exists a singularity fixing regular
homotopy connecting $f$ and $g$ that is fixed on $S^n_-$, i.e., $f
\sim_s g$ (rel $S^n_-$).
\end{thm}

We will use the following lemma in the proof of Theorem
\ref{thm:6}:

\begin{lem} \label{lem:13}
Let $f,g \colon S^n \to \Real^{2n-1}$ be locally generic maps such
that $S(f) = S(g) \subset \text{int}(S^n_+)$ and $f|S^n_- = g
|S^n_-$ are immersions. If $f \sim_s g$ then $f \sim_s g$ (rel
$S^n_-$).
\end{lem}

\begin{proof}
Choose a point $x \in \partial S^n_-$ and denote the differential
$(df)_x = (dg)_x$ by $j$. Let $S = S(f) = S(g)$ and $$\mathcal{B}
= \{\,h \in \imm(S^n_-, \Real^{2n-1}) \colon  (dh)_x = j\,\},$$
moreover $$\mathcal{E} = \{\,h \in L(S^n, \Real^{2n-1}) \colon
S(h) = S \,\,\text{and}\,\, (dh)_x = j \,\}.$$ We define a map $p
\colon \mathcal{E} \to \mathcal{B}$ by the formula $p(h) =
h|S^n_-$ for $h \in \mathcal{E}$.

We will now show that $p$ has the covering homotopy property. For
this end choose a subset $B \subset S^n_+$ such that $B$ is
diffeomorphic to an open $n$-disc centered at the North pole $N$
of $S^n$, moreover $S \subset B$. Denote by $R$ the intersection
of $S^n_+ \setminus B$ with the geodesic arc connecting $N$ and
$x$. Let $K = S^n_+ \setminus (B \cup R)$, then $\text{int}(K)$ is
diffeomorphic to an open $n$-disc. Let $P$ be a polyhedron, $\psi
\colon P \times I \to \mathcal{B}$ and $\phi_0 \colon P \times
\{0\} \to \mathcal{E}$, such that $\psi|(P \times \{0\}) = p \circ
\phi_0$. We must extend $\phi_0$ to $\phi \colon P \times I \to
\mathcal{E}$ so that $\psi = p \circ \phi$, i.e., the following
diagram is commutative:
\[
\begin{CD}
P \times I @>\phi>>
\mathcal{E} \\
@V\text{id}VV @VpVV \\
P \times I @>\psi >> \mathcal{B}.
\end{CD}
\]
For $r \in P$ and $t \in I$ we define $\phi(r,t)|(\overline{B}
\cup R)$ to be equal to $\phi_0(r)|(\overline{B} \cup R)$.
Moreover, let $(d \phi(r,t))|R = (d \phi_0(r))|R$. This is
possible since $(\overline{B} \cup R ) \cap S^n_- = \{x\}$ and
$(d\psi(r,t))_x = j$ for every $r \in P$ and $t \in I$. The map
$\phi_0(r)|K$ is an immersion for every $r \in P$ and $\phi(r,t)$
is already defined (together with normal derivatives) on $\partial
K$ for every $r \in P$ and $t \in I$. Thus using Smale's lemma
(Theorem \ref{thm:9}) $\phi(r,t)$ can be extended to $K$ as an
immersion for every $r \in P$ and $t \in I$.

It is well known that the space $\mathcal{B}$ is contractible,
thus its homotopy groups are all trivial. Take the fiber of $p$
defined by $\mathcal{F} = p^{-1}(f|S^n_-)$, then by assumption
$f,g \in \mathcal{F}$. Denote by $i$ the inclusion of
$\mathcal{F}$ into $\mathcal{E}$. Let us look at the following
part of the homotopy exact sequence of the fibration $p \colon
\mathcal{E} \xrightarrow{\mathcal{F}} \mathcal{B}:$
$$\pi_1(\mathcal{B})\longrightarrow \pi_0(\mathcal{F}) \xrightarrow{i_*} \pi_0(\mathcal{E})
\xrightarrow{p_*} \pi_0(\mathcal{B}).$$ Since $\pi_1(\mathcal{B})
= 0$ and $\pi_0(\mathcal{B}) = 0$, the morphism $i_*$ is an
isomorphism. By assumption $f \sim_s g$, moreover, a singularity
fixing regular homotopy $H$ connecting $f$ and $g$ can be chosen
so that $(dH_t)_x = j$ for every $t \in [0,1]$. Thus $i(f)$ and
$i(g)$ lie in the same path-component of $\mathcal{E}$, i.e., they
represent the same element of $\pi_0(\mathcal{E})$. Using the fact
that $i_*$ is a monomorphism we get that $f$ and $g$ lie in the
same path-component of $\mathcal{F}$. This means that there exists
a singularity fixing regular homotopy connecting $f$ and $g$ that
is fixed on $S^n_-$.
\end{proof}

\subsection{Merging double curves}

The idea of the procedure of eliminating the second obstruction
(and achieving the coincidence of $f$ and $g$ -- after a regular
homotopy -- also on the last part of $M^n$, i.e., on $U$) is as
follows: In \cite{Ekholm2} Ekholm showed that the regular homotopy
class of an immersion $S^n \to \Real^{2n-1}$ is completely
determined by the behavior of the map in a small neighborhood of
the double curves if $n \ge 4$. We are going to reduce the problem
to this result of Ekholm. First we  make $f$ and $g$ coincide in a
neighborhood of their singular sets. Then using Theorem
\ref{thm:8} we will get to the case of the immersions of Ekholm's
theorem. To manipulate the double curves of $f$ and $g$ we will us
the following construction.

In the presence of Whitney umbrella singular points we can merge
the closed double curves with the double curves connecting the
singular points. The merging is done by using the following
construction motivated by the Whitney-trick (for a reference see
\cite{Milnor} or \cite{Whitney1}). (Let us recall that in the
Whitney trick one defines a standard model and then one shows that
it can be embedded into the manifold under consideration.)

\begin{proof}[The construction]
In our case the standard model $S$ is chosen to be the standard
model of the Whitney-trick (two arcs in $D^2$ that intersect in
two points) times the interval $[-1,1]$ (see the left side of
Figure \ref{fig:3}). Thus we have two surfaces $P_1$ and $P_2$ in
$D^2 \times [-1,1]$, both diffeomorphic to the square, that
intersect in two line segments $l_1$ and $l_2.$ The regular
homotopy in this model is the identity near $\partial(D^2 \times
[-1,1])$ and is the isotopy of the standard model of the Whitney
trick ($l_2$ is pushed above $l_1$) near $D^2 \times \{0\}$. As a
result of this deformation the double points are removed near $D^2
\times \{0\}$, and thus the double curves of the final arrangement
coincide with the connected sum $l_1 \# l_2$ (see the right side
of Figure \ref{fig:3}).

\begin{figure}[t]
\includegraphics{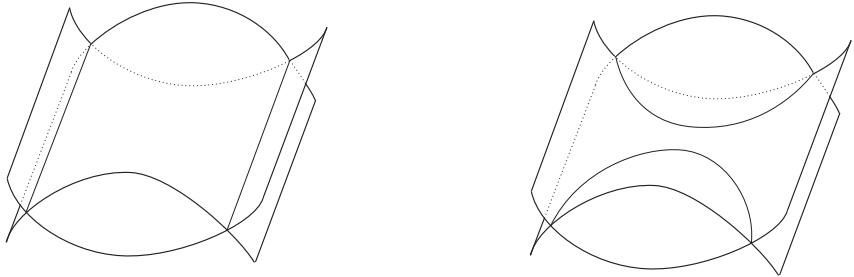}
\caption{The standard model} \label{fig:3}
\end{figure}

Let $p$ and $q$ be two double points of the generic map $f \colon
N^n \to \Real^{2n-1}$, where $N^n$ is any compact connected
manifold with boundary. All we have to do is to embed $D^2 \times
[-1,1] \times \Real^{n-2} \times \Real^{n-2}$ into $\Real^{2n-1}$
so that $P_1 \times \Real^{n-2} \times \{0\}$ and $P_2 \times
\{0\} \times \Real^{n-2}$ map to $\im(f)$, moreover $l_1 \times
\{0\}$ maps on the double curve through $p$ and $l_2 \times \{0\}$
maps on the double curve through $q$.

Let $f^{-1}(p) = \{p_1,p_2\}$ and $f^{-1}(q) = \{q_1,q_2\}$.
Choose an embedded curve $c_i$ in $N^n$ connecting $p_i$ and $q_i$
and let $C_i = f(c_i)$ for $i = 1,2$. We might suppose that $C_1
\cap C_2 = \{p,q\}$ and no other double or singular point of $f$
lies on $c_1$ and $c_2$. Thus $f|c_1$ and $f|c_2$ are embeddings.
Let $D_i$ be an embedded (tubular) neighborhood of $C_i$ in
$f(N^n)$ diffeomorphic to $D^n$ for $i=1,2$. The intersection $D_1
\cap D_2 = L_1 \sqcup L_2$, where $p \in L_1$ and $q \in L_2$.

Next we choose a Whitney disk $W$ containing $C_1$ and $C_2$ as in
\cite{Milnor}. Since $\dim(\im(f))+ \dim(W)= n + 2 < 2n-1$
(because $n \ge 4$), if $W$ is in general position then $W \cap
\im(f) = C_1 \cup C_2$. We might also suppose that $W$ is
transversal to $L_1$ and $L_2$. Since from now on we will work
only in a neighborhood of $W$, we might restrict our attention to
$D_1 \cup D_2$ and forget about the rest of $\im(f)$.

Construct a Riemannian metric on $\Real^{2n-1}$ for which $D_1
\perp D_2$. Choose an orientation of $D_1$ and a co-orientation of
$D_2$ in $\Real^{2n-1}$. Let $$\alpha_1, \dots, \alpha_{n-1} \in
\nu_p(D_2 \subset \Real^{2n-1}) \subset T_pD_1$$ be a positive
normal frame of $D_2$. Let $\tilde{\nu}(p) \in T_pL_1$ be taken so
that $\tilde{\nu}(p),\alpha_1, \dots, \\ \alpha_{n-1}$ is a
positive basis of $T_pD_1$. Similarly, let $$\beta_1, \dots,
\beta_{n-1} \in \nu_q(D_2 \subset \Real^{2n-1})\subset T_qD_1$$ be
a positive normal frame of $D_2$ and $\tilde{\nu}(q) \in
T_qL_2$ be chosen so that $\tilde{\nu}(q), \beta_1, \dots, \\
\beta_{n-1}$ is a \emph{negative} basis of $T_qD_1$.

We can extend $\tilde{\nu}$ to $C_1 \cup C_2$ so that
$\tilde{\nu}|C_i$ is tangent to $D_i$ since $n \ge 3$. Next we
extend $\tilde{\nu}$ from $C_1 \cup C_2$ to a normal field of $W$:
Since $\nu(W \subset \Real^{2n-1})$ is the trivial bundle and $C_1
\cup C_2$ is homotopy equivalent to $S^1$, the obstruction to
extend $\tilde{\nu}$ is an element of $\pi_1(V_{2n-3,1}) =
\pi_1(S^{2n-4})$. But this group is $0$ since $2n-4 > 1$ if $n \ge
3$. Taking the exponential map along $\tilde{\nu}$ on $W$ we can
embed the standard model $S = D^2 \times [-1,1]$ into
$\Real^{2n-1}$. Denote the image of $S$ by $U$.

Orient $C_1$ from $p$ to $q$ and denote by $\tau(x)$ the positive
unit tangent vector of $C_1$ at $x$. Moreover, let $\nu(p) =
\tau(p)$ and $\nu(q) = -\tau(q)$. Since $n \ge 4$ we can extend
$\nu$ along $C_2$ as a normal field of $D_2$.

Choose vectors $\xi_1(p), \dots, \xi_{n-2}(p)$ in $T_pD_1$ so that
$\tilde{\nu}(p), \nu(p), \xi_1(p), \dots, \xi_{n-2}(p)$ is a
positive basis of $T_pD_1$. Using parallel translation extend
$\xi_1, \dots, \xi_{n-2}$ to $C_1$. Then $$\tilde{\nu}(q), \nu(q),
\xi_1(q),\dots, \xi_{n-2}(q) \in T_qD_1 $$ is a negative basis
since $\nu(q) = -\tau(q)$. Thus from the definition of
$\tilde{\nu}(q)$ we get that $\nu(q),\xi_1(q), \dots,
\xi_{n-2}(q)$ is a positive basis of $\nu_q(D_2 \subset
\Real^{2n-1})$. Moreover,  $\nu(p),\xi_1(p), \dots, \xi_{n-2}(p)$
is a positive basis of $\nu_p(D_2 \subset \Real^{2n-1})$ and
$\nu|C_2$ is a continuous vector field. All this together implies
that the $(n-2)$-frames $\underline{\xi}(p)$ and
$\underline{\xi}(q)$ can be extended onto $C_2$ to be transversal
to $D_2$ and $\nu$.

We can extend the frame $\underline{\xi}$ to $U$: The pair $(U,
C_1\cup C_2)$ is homotopy equivalent (deformation retracts onto)
$(D^2,S^1)$. Thus the obstruction lies in $\pi_1(V_{2n-4,n-2}) =
0$ since $n \ge 4$ implies that $1 < (2n-4) - (n-2)$.

Finally we construct an $(n-2)$-frame $\underline{\eta}$ on $U$
that is tangent to $D_2$ on $C_2$. Look at the bundle of
orthogonal $(n-2)$-frames on $U$ that are perpendicular to $U$ and
$\underline{\xi}$. Since $U$ is contractible, this bundle is
trivial; let $\underline{\eta}$ be any section. Notice that $\nu(U
\cap D_2 \subset U)$ is spanned by the vector field $\nu|C_2$ (the
subset $U \cap D_2$ corresponds to $P_2$ in the standard model
$S$). Because $\underline{\eta}|C_2 \perp \underline{\xi}, U$, we
get that $$\underline{\eta}|C_2 \perp \langle \underline{\xi},
\nu|C_2 \rangle = \nu(D_2 \subset \Real^{2n-1})|C_2.$$ Thus
$\underline{\eta}|C_2$ is tangent to $D_2$, as required.

Using the exponential map on the fields $\underline{\xi}$ and
$\underline{\eta}$ along $U$ we can embed $S \times \Real^{n-2}
\times \Real^{n-2}$ into $\Real^{2n-1}$ as described above.
\end{proof}

The regular homotopy constructed above is clearly singularity
fixing.

\subsection{The proof of Theorem \ref{thm:6}}

Using Lemma \ref{lem:13} we only have to prove that $f \sim_s g$.
With a singularity fixing regular homotopy we can perturb $f$ and
$g$ to obtain generic maps. Let us recall (see the beginning of
the previous section) that the main idea of the proof is to reduce
the problem to the case of immersions and then apply Ekholm's
theorem.

Let $l_f$ (respectively $l_g$) denote the closure of the double
curve of $f$ (respectively $g$) that connects the two singular
points of $f$ (respectively $g$). With the aid of the above
construction we merge using a singularity fixing regular homotopy
all closed double curves of $f$ (respectively $g$) with $l_f$
(respectively $l_g$). From now on we will suppose that $l_f$ is
the only double curve of $f$ and $l_g$ is the only double curve of
$g$. Let $m_f = f^{-1}(l_f)$ and $m_g = g^{-1}(l_g)$. Since $n \ge
4$ there exists a diffeotopy $\{\,d_t \colon t \in [0,1]\,\}$ of
$S^n$ that fixes $S(f)=S(g)$ and takes $m_g$ to $m_f$. Composing
$f$ with $\{d_t\}$ we can suppose that $m_f = m_g$. Let $m = m_f =
m_g$ and $S(f)=S(g)= \{p,q\}$. If $n$ is even we have to be
careful when choosing the diffeotopy $\{d_t\}$: If we order the
two components of $m \setminus \{p,q\}$ then this defines signs of
the singular points of $f$ and $g$ (see Definition \ref{defn:5}),
and it might happen that $S_+(f) \neq S_+(g)$ (the equality
$\#S(f) = \#S(g) = 0$ still holds). In this case we compose $f$
with a diffeotopy of $S^n$ that fixes $p$ and $q$ and swaps the
two components of $m \setminus \{p,q\}$ in order to swap $S_+(f)$
and $S_-(f)$. This is possible if $n \ge 4$. Thus if $n$ is even
we also fix an ordering of the components of $m \setminus \{p,q\}$
and suppose that for the induced signs $S_{\pm}(f) = S_{\pm}(g)$.

Next we  compose $f$ with a diffeotopy of $\Real^{2n-1}$ that
moves $l_f$ to $l_g$. Thus we can suppose that $f|m = g|m$ and let
$l = l_f = l_g$. Using Lemma 2.14 of \cite{Ekholm} we can
straighten $f$ and $g$ close to their self intersection
$\text{int}(l)$, i.e., we can achieve that they agree with their
normal derivatives on a neighborhood of $m$.

We now sketch the idea of the construction that follows: The
singular points of $f$ and $g$ coincide and have the same signs if
$n$ is odd. The sign is the only local isotopy invariant a
Whitney-umbrella singularity can have. Moreover, any two bundles
over an interval are equivalent. Since $m$ is an orientation
preserving curve in $S^n$ the two cross-caps at the end of $l_f$
and $l_g$ are in the same "relative position". So $f$ and $g$ are
isotopic if restricted to a small neighborhood of $m$, fixing
$S(f) = S(g)$. This isotopy can be extended to an ambient isotopy
of $\Real^{2n-1}$. So we might suppose that $f$ and $g$ agree in a
small open tubular neighborhood $V$ of their common double locus
$m$. The exact details are as follows:

Since $S^n$ is orientable, the normal bundle $\nu(m \subset S^n)$
is trivial. Thus we can choose a normal framing $\underline{\nu}$
of $m$ in $S^n$. Notice that $\ker(df_p)= \ker(dg_p) = T_pm$ and
$\ker(df_q)= \ker(dg_q) = T_qm$, so $df(\underline{\nu}_p)$ and
$dg(\underline{\nu}_p)$ are $(n-1)$-frames in
$T_{f(p)}\Real^{2n-1} = T_{g(p)}\Real^{2n-1}$. Similarly,
$df(\underline{\nu}_q)$ and $dg(\underline{\nu}_q)$ are
$(n-1)$-frames in $T_{f(q)}\Real^{2n-1} = T_{g(q)}\Real^{2n-1}$.
Moreover, $\underline{\nu}_f = df(\underline{\nu}|(m \setminus
\{p,q\}))$ and $\underline{\nu}_g = dg(\underline{\nu}|(m
\setminus \{p,q\}))$ provide $(2n-2)$-frames along
$\text{int}(l)$.

Suppose that $S_+(f) = S_+(g) = \{p\}$ and $S_-(f) = S_-(g) =
\{q\}$. Let $a = f(p) = g(p)$ and $b = f(q) = g(q)$, then $a$ and
$b$ are the endpoints of the curve $l$. Orient $l$ from $a$ to
$b$, this induces an orientation of $\nu(l \subset \Real^{2n-1})$.
(Recall that we fixed an orientation of $\Real^{2n-1}$ to define
the signs of the Whitney-umbrellas.) Then from Definition
\ref{defn:5} of the signs of the singular points we can see that
$\left\{\,\underline{\nu}_f, \underline{\nu}_g\,\right\}$ provides
a negative basis of $\mu = \nu\left(\text{int}(l) \subset
\Real^{2n-1}\right)$. Thus there exists a homotopy $T \colon
\text{int}(l) \times I \to GL_+(2n-2, \mathbb{R})$ such that for
every $x \in \text{int}(l)$ the transformation $T(x,0) =
\text{id}_{\Real^{2n-2}}$ and $(\underline{\nu}_f)_x T(x,1) =
(\underline{\nu}_g)_x$. The homotopy $T$ can be extended to an
ambient isotopy of $f$ that takes $f$ to $g$ on a  small open
tubular neighborhood $V$ of $m$, and this isotopy is identical
outside a neighborhood of $l$. Thus we can suppose that $f|V =
g|V$.

Let $D \subset V$ be diffeomorphic to the closed $n$-disc $D^n$,
so that $p,q \in \text{int}(D)$ and $f|\partial D$ ($= g |\partial
D$) is an embedding. Since the only double locus of $f$ and $g$ is
$m \subset V$, the distances $d_1 = d(f(D), f(S^n \setminus V))$
and $d_2 = d(g(D), g(S^n \setminus V))$ are positive. Let $0 <
\varepsilon < \min(d_1,d_2)$. Using Theorem \ref{thm:8} choose a
generic immersion $h \colon D \to \Real^{2n-1}$ that agrees with
$f|\partial D$ in a neighborhood of $\partial D$, moreover the
$C^0$-distance of $h$ and $f|D$ is $< \varepsilon$. (This is
possible since $\#S(f|D) = 0$ and $f|\partial D$ is an embedding.)
Denote by $f_1$ (respectively $g_1$) the immersion of $S^n$ in
$\Real^{2n-1}$ that agrees with $f$ (respectively $g$) on $S^n
\setminus D$ and with $h$ on $D$. From the choice of $h$ we can
see that $f(S^n \setminus V) \cap h(D) = \emptyset$ and $g(S^n
\setminus V) \cap h(D) = \emptyset$. Moreover, since $f|(S^n
\setminus V)$ and $g|(S^n \setminus V)$ are embeddings that do not
intersect $f(V \setminus D) = g(V \setminus D)$, the double loci
of $f_1$ and $g_1$ are contained in $V$. But $f_1|V = g_1|V$, thus
$f_1$ and $g_1$ agree in a neighborhood of their common double
locus. Theorem 1 in \cite{Ekholm2} implies that the regular
homotopy class of a generic immersion $S^n \to \Real^{2n-1}$ can
be expressed in terms of the geometry of the self intersection if
$n \ge 4$ . Thus $f_1 \sim g_1$.

Since $f_1|D = g_1|D$, the regular homotopy connecting $f_1$ and
$g_1$ can be chosen to be constant on $D$, i.e., $f_1 \sim g_1$
(rel $D$). This gives us a regular homotopy $f \sim g$ (rel $D$),
thus $f \sim_s g$.

\subsection{An application}

In part 4 of \cite{Juhasz} the $n=2$ case of Theorem \ref{thm:1}
was applied to projections of immersions of surfaces into
$\Real^4$. Knowing Theorem \ref{thm:1} the same results can be
generalized without modification of the proof.

Let us first recall some definitions.

\begin{defn}
Let
\[
\imm_{\pi}(M^n, \Real^{2n}) = \left\{\, F \in \imm(M^n,
\Real^{2n}) \colon \pi \circ F \in L(M^n,
\Real^{2n-1})\,\right\}
\]
be the subspace of $\imm(M^n, \Real^{2n})$ formed by those
immersions whose projection in $\Real^{2n-1}$ is locally generic.
Two immersions $F,G \in \imm_{\pi}(M^n, \mathbb{R}^{2n})$ are
called \emph{$\pi$-homotopic} (denoted by $F \sim_{\pi} G$) if
they are in the same path-component of $\imm_{\pi}(M^n,
\mathbb{R}^{2n})$.
\end{defn}

The generalization of Theorem 4.7 of \cite{Juhasz} is then the
following.

\begin{thm}
If $n=2$ or $n>3$ and $F, G \in \imm_{\pi}(M^n, \mathbb{R}^{2n})$
then
\[ F \sim_{\pi} G \Leftrightarrow [F \sim G \,\, \text{and} \,\,
\pi \circ F \sim \pi \circ G].
\]
\end{thm}

\section*{Acknowledgement}

I would like to take this opportunity to thank Professor Tobias
Ekholm for valuable advice and Professor Andr\'as Sz\H{u}cs for
reading earlier versions of this paper and making several useful
remarks.

% ----------------------------------------------------------------
\bibliographystyle{amsplain}
\bibliography{topology}
\end{document}